\documentclass[aop,MSNbibl,dvips]{arximspdf}
\usepackage{graphicx}
\doi{10.1214/12-AOP755} 
\volume{42}
\issue{1}
\pubyear{2014}
\firstpage{1}
\lastpage{40}

\makeatletter
\newcommand{\rrvert}{\vert}
\newcommand{\llvert}{\vert}

\newtheorem{Theorem}{Theorem}[section]
\newproclaim{Definition}{Definition}[section]
\newtheorem{Lemma}[Theorem]{Lemma}
\newtheorem{Corollary}[Theorem]{Corollary}
\newtheorem{Proposition}[Theorem]{Proposition}
\newproclaim{Remark}[Theorem]{Remark}

\newcommand{\cal}{\mathcal}
\newcommand{\I}{{\mathbf{I}}}
\newcommand{\FF}{{\mathcal F}}
\newcommand{\prob}{{\mathbf P}}

\def\eps{\varepsilon}
\def\G{{\cal{G}}}
\def\phi{\varphi}
\def\eps{\varepsilon}
\def\inj{\operatorname{inj}}
\newcommand{\U}{{\mathcal U}}

\makeatother

\begin{document}
\begin{frontmatter}

\title{Asymptotic behavior and distributional limits of preferential attachment graphs}
\runtitle{Preferential attachment limits}

\begin{aug}
\author[A]{\fnms{Noam} \snm{Berger}\ead[label=e1]{berger@math.huji.ac.il}},
\author[B]{\fnms{Christian} \snm{Borgs}\ead[label=e2]{borgs@microsoft.com}},
\author[B]{\fnms{Jennifer T.} \snm{Chayes}\ead[label=e3]{chayes@microsoft.com}}\\
\and
\author[C]{\fnms{Amin} \snm{Saberi}\corref{}\ead[label=e4]{saberi@stanford.edu}}
\runauthor{Berger, Borgs, Chayes and Saberi}
\affiliation{Hebrew University, Microsoft Research, Microsoft
Research and Stanford~University}
\address[A]{N. Berger\\
Mathematics Department\\
Hebrew University\\
Jerusalem 91904\\
Israel\\
\printead{e1}}
\address[B]{C. Borgs\\
J. T. Chayes\\
Microsoft Research New England\\
Cambridge, Massachusetts 02142\\
USA\\
\printead{e2}\\
\hphantom{E-mail: }\printead*{e3}}
\address[C]{A. Saberi\\
Management Science and Engineering\\
Stanford University\\
Palo Alto, California 94305\\
USA\\
\printead{e4}}
\end{aug}

\received{\smonth{9} \syear{2010}}
\revised{\smonth{3} \syear{2012}}

%
\begin{abstract}
We give an explicit construction of the weak local limit
of a class of preferential attachment graphs. This limit contains all
local information and allows several computations that are otherwise
hard, for example, joint degree distributions and, more generally, the
limiting distribution of subgraphs in balls of any given radius $k$
around a random vertex in the preferential attachment graph. We also
establish the finite-volume corrections which give the approach to the
limit.
\end{abstract}

%
\begin{keyword}[class=AMS]
\kwd{60C05}
\kwd{60K99}
\end{keyword}
\begin{keyword}
\kwd{Preferential attachment graphs}
\kwd{graph limits}
\kwd{weak local limit}
\end{keyword}

\end{frontmatter}

\section{Introduction}

About a decade ago, it was realized that the Internet has a
power-law degree distribution \cite{Faloutsos,AJB99}. This
observation led to the so-called preferential
attachment model of Barab\'asi and Albert \cite{Barabasi1},
which was later used to explain the observed power-law degree
sequence of a host of real-world networks, including social and
biological networks, in addition to technological ones.
The first rigorous analysis of a preferential
attachment model, in particular proving that it has small
diameter, was
given by Bollob\'as and Riordan \cite{brdiam}.
Since these works there has been a tremendous amount of study,
both nonrigorous and rigorous, of the random graph models that
explain the power-law degree distribution; see \cite{BAreview}
and \cite{BRreview} and references therein for some of the
nonrigorous and rigorous work, respectively.

Also motivated by the growing graphs appearing in real-world
networks, for the past five years or so, there has been much
study in the mathematics community of notions of graph limits.
In this context, most of the work has focused on dense graphs.
In particular, there has been a series of papers on a notion
of graph limits defined via graph homomorphisms
\cite{BCLSV-rev,dense1,dense2,LSz}; these have been shown to be
equivalent to limits defined in many other senses
\cite{dense1,dense2}. Although most of the results in this work
concern dense graphs, the paper \cite{BCLSV-rev} also
introduces a notion of graph limits for sparse graphs with
bounded degrees in terms of graph homomorphisms; using
expansion methods from mathematical physics, Borgs et al. \cite{sparse}
establishes some general results on this type of limit for
sparse graphs. Another recent work \cite{BR07} concerns limits
for graphs which are neither dense nor sparse in the above
senses; they have average degrees which tend to infinity.

Earlier, a notion of a weak local limit of a sequence of graphs
with bounded degrees was given by Benjamini and Schramm
\cite{BS01} (this notion was in fact already implicit in
\cite{Ald99}). Interestingly, it is not hard to show that the
Benjamini--Schramm limit coincides with the limit defined via
graph homomorphisms in the case of sparse graphs of bounded
degree; see \cite{Elek} for yet another equivalent notion of
convergent sequences of graphs with bounded degrees.

As observed by Lyons \cite{Lyo05}, the notion of graph
convergence introduced by Benjamini and Schramm is meaningful
even when the degrees are unbounded, provided the \textit{average
degree} stays bounded. Since the average degree of the
Barab\'asi--Albert graph is bounded by construction, it is
therefore natural to ask whether this graph sequence converges
in the sense of Benjamini and Schramm.

In this paper, we establish the existence of the
Benjamini--Schramm limit for the Barab\'asi--Albert graph by
giving an explicit construction of the limit process, and use
it to derive various properties of the
limit.
Our results
cover the case of both uniform and preferential attachments
graphs.\footnote{Note, however, that we do not cover models
exhibiting densification in the sense of Leskovec, Kleinberg and Faloutsos
\cite{LKF07}; see \cite{LCKFG10} for
a mathematical model exhibiting this phenomenon.
Indeed, these models are outside the scope of convergence considered in
this paper,
since they have bounded diameter and growing average degree,
and hence do not converge in the
sense of Benjamini--Schramm.}
Moreover, our methods establish the finite-volume
corrections which give the approach to the limit.

Our proof uses a representation, which we first introduced in
\cite{BBCS05}, to analyze processes that model the spread of
viral infections on preferential attachment graphs. Our
representation expresses the preferential attachment model
process as a combination of several P\'olya urn processes. The
classic P\'olya urn model was of course proposed and analyzed in
the beautiful work of P\'olya and Eggenberger in the early
twentieth century \cite{EP}; see
\cite{durrett} for a basic reference.
Despite the fact that our P\'olya urn representation is a
priori only valid for a limited class of preferential
attachment graphs, we give an approximating coupling which
proves that the limit constructed here is the limit of a much
wider class of preferential attachment graphs.

Our alternative representation contains much more independence
than previous representations of preferential attachment and is
therefore simpler to analyze. In order to demonstrate this, we
also give a few applications of the limit. In particular, we
use the limit to calculate the degree distribution and the
joint degree distribution of a typical vertex with
the vertex it attached to in the preferential attachment
process (more precisely, a vertex chosen uniformly from the
ones it attached to).

\section{Definition of the model and statements of results}\label{secresults}

\subsection{Definition of the model}
\label{secdef-mod}

The preferential attachment graph we define generalizes the
model introduced by Barab\'asi and Albert \cite{Barabasi1} and
rigorously constructed in \cite{brdiam}. Fix an integer $m\geq
2$ and a real number $0\leq\alpha<1$. We will construct a
sequence of graphs $(G_n)$ (where $G_n$ has $n$ vertices
labeled $1,\ldots,n$) as follows:

$G_1$ contains one vertex and no edges, and $G_2$ contains two
vertices and $m$ edges connecting them. Given $G_{n-1}$ we
create $G_n$ the following way: we add the vertex $n$ to the
graph, and choose $m$ vertices $w_1,\ldots,w_m$, possibly with
repetitions, from $G_{n-1}$. Then we draw edges between $n$ and
each of $w_1,\ldots,w_m$. Repetitions in the sequence
$w_1,\ldots,w_m$ result in multiple edges in the graph $G_n$.

We suggest three different ways of choosing the
vertices $w_1,\ldots,w_m$. The first two ways, the independent and
the conditional, are natural ways which we consider of
interest, and are the two most common interpretations of the
preferential attachment model. The third way, that is, the
sequential model, is less natural, but is much easier to
analyze because it is exchangeable, and therefore by
de-Finetti's theorem (see \cite{durrett}) has an alternative
representation, which contains much more independence. We call
this representation the P\'olya urn representation because the
exchangeable system we use is the P\'olya urn scheme.
\begin{enumerate}[(2)]
\item[(1)]\label{itemindep} The independent model:
$w_1,\ldots,w_m$ are chosen independently of each other
conditioned on the past, where for each $i=1,\ldots,m$,
we choose $w_i$ as follows: with probability
$\alpha$, we choose $w_i$ uniformly from the vertices
of $G_{n-1}$, and with probability $1-\alpha$, we
choose $w_i$ according to the preferential attachment
rule, that is, for all $k=1,\ldots, n-1$,
\[
\prob(w_i=k )=\frac{\deg_{n-1}(k)}Z,
\]
where $Z$ is the normalizing constant
$Z=\sum_{k=1}^{n-1}\deg_{n-1}(k)= 2m(n-2)$.

\item[(2)]\label{itemindepcond} The conditional model: here we
start with some predetermined graph structure for the
first $m$ vertices. Then at each step, $w_1,\ldots,w_m$
are chosen as in the independent case, \emph{conditioned} on them being
different from one another.

\item[(3)]\label{itemsequential} The sequential model:
$w_1,\ldots,w_m$ are chosen inductively as follows: with
probability
$\alpha$, $w_1$ is chosen uniformly, and with probability
$1-\alpha$, $w_1$ is chosen according to the preferential
attachment rule, that is, for every $k=1,\ldots,n-1$, we take
$w_1=k$ with probability $(\deg_{n-1}(k))/Z$ where as
before $Z=2m(n-2)$. Then we proceed inductively, applying
the same rule, but with two modifications:
\begin{enumerate}[(a)]
\item[(a)] When determining $w_i$, instead of the degree
$\deg_{n-1}(k)$, we use
\[
\deg^\prime_{n-1}(k)=\deg_{n-1}(k)+\#\{1\leq j\leq
i-1 \mid w_j=k\}
\]
and normalization constant
\[
Z^\prime=\sum_{k=1}^{n-1}\bigl(
\deg^\prime_{n-1}(k)\bigr)= 2m(n-2)+i-1.
\]

\item[(b)] The probability of uniform connection will be
%
%
\begin{equation}
\label{tilde-alpha} \tilde\alpha=\alpha\frac{2m(n-1)}{2m(n-2)+2m\alpha
+(1-\alpha)(i-1)} =\alpha+O
\bigl(n^{-1}\bigr)
\end{equation}
%
rather than $\alpha$.
\end{enumerate}
\end{enumerate}

We will refer to all three models as versions of the
preferential attachment graph, or PA-graph, for short. Even
though we consider the graph $G_n$ as undirected, it will often
be useful to think of the vertices $w_1,\ldots,w_m$ as vertices
which ``received an edge'' from the vertex $n$, and of $n$ as a
vertex which ``sent out $m$ edges'' to the vertices
$w_1,\ldots,w_m$. Note in particular, that the degree of a
general vertex $v$ in $G_n$ can be written as $m+q$, where $m$
is the number of edges sent out by $v$ and $q$ is the
(random) number of edges received by $v$.

\subsection{P\'olya urn representation of the sequential model}

Our first theorem gives the P\'olya urn representation of the
sequential model. To state it, we use the standard notation
$X\sim\beta(a,b)$ for a random variable $X\in[0,1]$ whose
density is equal to $\frac1Z x^{a-1}(1-x)^{b-1}$, where
$Z=\int_0^1 x^{a-1}(1-x)^{b-1}\,dx$. We set
\[
u=\frac\alpha{1-\alpha}.
\]
Note that $u\in[0,\infty)$.

\begin{Theorem}
\label{thm1} Fix $m$, $\alpha$ and $n$. Let $\psi_1=1$, let
$\psi_2,\ldots,\psi_n$ be independent random variables with
%
%
\begin{equation}
\label{psik-dis} \psi_j\sim\beta\bigl(m+2mu, (2j-3)m+2mu (j-1)
\bigr)
\end{equation}
and let
%
%
\begin{equation}
\label{Sk} \phi_j=\psi_j\prod
_{i=j+1}^n(1-\psi_i),\qquad S_k=
\sum_{j=1}^k\phi_j\quad
\mbox{and}\quad I_k=[S_{k-1},S_k).
\end{equation}
Conditioned on $\psi_1,\ldots,\psi_n$,
choose $\{U_{k,i}\}_{k=1,\ldots, n, i=1,\ldots,m}$ as a sequence
of independent random variables, with $U_{{ k,i}}$ chosen
uniformly at random from $[0,S_{k-1}]$. Join two vertices $j$
and $k$ if $j<k$ and $U_{k,i}\in I_j$ for some
$i\in\{1,\ldots,m\}$ (with multiple edges between $j$ and $k$ if
there are several such $i$). Denote the resulting random
multi-graph by $G_n$.

Then $G_n$ has the same distribution as the sequential
PA-graph.
\end{Theorem}

\begin{figure}

\includegraphics{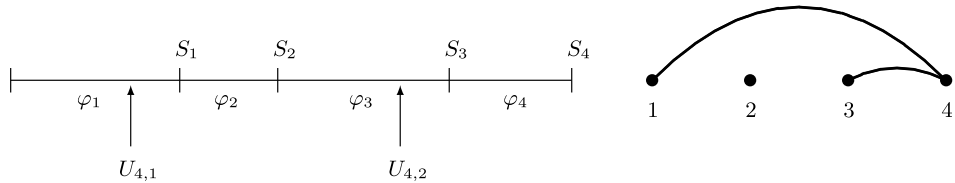}

\caption{The P\'olya-representation of
the sequential model for $m=2$, $n=4$ and $k=4$. The variables
$U_{4,1}$ and $U_{4,2}$ are chosen uniformly at random from
$[0,S_3]$.}\label{fig1}
\end{figure}

Figure \ref{fig1} illustrates this theorem.

%

It should be noted that the $\alpha= 0$ case of the sequential
model defined here differs slightly from the model of
Bollob\'as and Riordan \cite{brdiam} in that they allow
(self-)loops, while we do not. In fact, a minor alteration of
our P\'olya urn representation models their graph, and we
suspect that a minor alteration of their pairing representation
can model our graph.

\subsection{Definition of the P\'olya-point graph model}

\subsubsection{Motivation}
\label{secexplore}
%

The Benjamini--Schramm notion \cite{BS01} of weak convergence
involves the view of the graph $G_n$ from the point of view of
a ``root'' $k_0$ chosen uniformly at random from all vertices
in $G_n$. More precisely, it involves the limit of the
sequence of balls of radius $1,2,\ldots\,$, about this root; see
Definition \ref{defBS-limit} in Section \ref{secmainresult}
below for the details.

It turns out that for the sequential model, this limit is
naturally described in terms of the random variables $S_{k-1}$
introduced in Theorem \ref{thm1}. To explain this, it is
instructive to first consider the ball of radius $1$ around the
random root $k_0$. This ball contains the $m$ neighbors of
$k_0$ that were born before $k_0$ and received an edge from
$k_0$ under the preferential attachment rule described above,
as well as a random number $q_0$ of neighbors that were born
after $k_0$ and send an edge to $k_0$ at the time they were
born. We denote these neighbors by $k_{01},\ldots,k_{0m}$ and
$k_{0,m+1},\ldots,k_{0,m+q_0}$, respectively.

Let
%
%
\begin{equation}
\label{eqdefu} \chi=\frac{1+2u}{2+2u} \quad\mbox{and}\quad \psi=\frac{1-\chi
}{\chi}{ =
\frac1{1+2u}}
\end{equation}
and note that $\frac12\leq\chi<1$ and $0<\psi\leq1$. As we will
see, the random variables $S_{k-1}$ behave asymptotically like
$(k/n)^\chi$, implying in particular that the distribution of
$S_{k_0-1}$ tends to that of a random variable $x_0=y_0^\chi$,
where $y_0$ is chosen uniformly at random in $[0,1]$. The
limiting distribution of $S_{k_{01}-1},\ldots,S_{k_{0m}-1}$
turns out to be quite simple as well: in the limit these random
variables are i.i.d. random variables $x_{0,i}$ chosen uniformly
from $[0,x_0]$, a distribution which is more or less directly
inherited from the uniform random variables $U_{k,i}\in
[0,S_{k_0-1}]$ from Theorem \ref{thm1}. The limiting
distribution of the random variables
$S_{k_{0,m+1}-1},\ldots,S_{k_{0,m+q_0}-1}$ is slightly more
complicated to derive and is given by a Poisson process in
$[x_0,1]$ with intensity
\[
{ \gamma_0\frac{\psi x^{\psi-1}}{x_0^{\psi}}\,dx.}
\]
Here $\gamma_0$ is a random ``strength'' which arises as a limit
of the $\beta$-distributed random variable $\psi_{k_0}$, and is
distributed according to $\Gamma({m+2mu},1)$. Here, as usual,
$\Gamma(a,b)$ is used to denote a distribution on $[0,\infty)$
which has density $\frac1{Z} x^{a-1}e^{-bx}$, with
$Z=\int_{0}^\infty x^{a-1}e^{-bx}\,dx$.

Next, consider the branching that results from exploring the
neighborhood of a random vertex in $G_n$ in a ball of radius
bigger than one. In each step of this exploration, we will find
two kinds of children of the current vertex $k$: those which
were born before $k$, and were attached to $k$ at the birth of
$k$, and those which were born after $k$, and were connected to
$k$ at their own birth. There are always either $m$ or $m-1$
children of the first kind (if $k$ was born after its parent,
there will be $m-1$ such children, since one of the $m$ edges
sent out by $k$ was sent out to $k$'s parent; otherwise there
will be $m$ children of the first type). The number of children
of the second kind is a random variable.

In the limit $n\to\infty$, this branching process leads to a
random tree whose vertices, $\bar a$, carry three labels: a
``strength'' $\gamma_{\bar a}\in(0,\infty)$ inherited from the
$\beta$-random variables $\psi_k$, a ``position'' $x_{\bar
a}\in[0,1]$ inherited from the random variables $S_{k-1}$ and
a type which can be either { $L$ (for ``left'') or $R$ (for
``right'')}, reflecting whether the vertex $k$ was born before or
after its parent. While the strengths of vertices of type { $R$}
turn out to be again $\Gamma({m+2mu},1)$-distributed, this is
not the case for vertices of type { $L$}, since a vertex with
higher values of $\psi_k$ has a larger probability of receiving
an edge from its
child.
In the limit, this will be reflected
by the fact that the strength of vertices of type { $L$} is
$\Gamma({m+2mu}+1,1)$-distributed.

\subsubsection{Formal definition}
\label{secpolyapointdef}


The main goal of the previous subsection was to give an
intuition of the structure of the neighborhood of a random
vertex. We will show that asymptotically, the branching
process obtained by exploring the neighborhood of a random
vertex $k_0$ in $G_n$ is given by a random tree with a certain
distribution. In order to state our main theorem, we give a
formal definition of this tree.

Let $F$ be the Gamma distribution $\Gamma({ m+2mu},1)$, and let
$F^\prime$ be the Gamma distribution $\Gamma({ m+2mu+1},1)$. We
define a random, rooted tree $(T,0)$ with vertices labeled by
finite sequences
\[
\bar{a}=(0,a_1,a_2,\ldots,a_l)
\]
inductively as follows:
\begin{itemize}
\item The root $(0)$ has a position $x_0=y_0^\chi$, where
$y_0$ is chosen uniformly at random in $[0,1]$. In the rest of the paper,
for notational convenience, we will write $0$ instead of $(0)$ for the root.

\item In the induction step, we assume that
$\bar{a}=(0,a_1,a_2,\ldots,a_l)$ and the corresponding
variable $x_{\bar{a}}\in[0,1]$ have been chosen in a
previous step. Define $(\bar a,j)$ as
$(0,a_1,a_2,\ldots,a_l,j)$, $j=1,2,\ldots\,$, and set
\[
m_-({\bar{a}})= \cases{ %
m, &\quad if $\bar a$ is the root or of type $L$,
\cr
m-1, &\quad if $\bar a$ is of type $R$.}
\]
We then take
\[
\gamma_{\bar{a}}\sim\cases{F, &\quad if $\bar a$ is the root or of type $R$,
\cr
F', &\quad if $\bar a$ is of type $L$,}
\]
independently of everything previously sampled, choose
$x_{(\bar{a},1)},\ldots,x_{(\bar{a},m_-(\bar{a}))}$ i.i.d.
uniformly at random in $[0,x_{\bar{a}}]$, and
$x_{(\bar{a},m_-(\bar{a})+1)},\ldots,x_{(\bar{a},m_-(\bar{a})+q
_{\bar{a}})}$
as the points of a Poisson process with intensity
%
%
\begin{equation}
\label{poisson-intensity} \rho_{\bar a}(x) \,dx = 
{
\gamma_{\bar{a}}\frac{\psi x^{\psi-1}}{x_{\bar{a}}^{\psi}}\,dx}
\end{equation}
on $[x_{\bar{a}},1]$ (recall that $0<\psi\leq1$).
The children of $\bar{a}$ are the vertices
$(\bar{a},1),\ldots,\break
(\bar{a},m_-(\bar{a})+q_{\bar{a}})$, with
$(\bar{a},1),\ldots,(\bar{a},m_-(\bar{a}))$ called of
type $L$, and the remaining ones called of type $R$.
\end{itemize}

We continue this process ad infinitum to obtain an infinite,
rooted tree $(T,0)$. We call this tree the P\'olya-point graph, and
the point process $\{x_{\bar a}\}$ the P\'olya-point process.

\subsection{Main result}
\label{secmainresult}

We are now ready to formulate our main result, which states
that in all three versions, the graph $G_n$ converges to the
P\'olya-point graph in the sense of \cite{BS01}.

Let $\G$ be the set of rooted graphs, that is, the set of all
pairs $(G,x)$ consisting of a connected graph $G$ and a
designated vertex $x$ in $G$, called the root. Two rooted
graphs $(G,x),(G',x')\in\G$ are called isomorphic if there is
an isomorphism from $G$ to $G'$ which maps $x$ to $x'$. Given
a finite integer $r$, we denote the rooted ball of radius $r$
around $x$ in $(G,x)\in\G$ by $B_r(G,x)$. We then equip $\G$
with the $\sigma$-algebra generated by the events that
$B_r(G,x)$ is isomorphic to a finite, rooted graph $(H,y)$
(with $r$ running over all finite, positive integers, and
$(H,y)$ running over all finite, rooted graphs),\vadjust{\goodbreak} and call
$(G,x)$ a random, rooted graph if it is a sample from a
probability distribution on $\G$.
We write $(G,x)\sim(G',x')$ if $(G,x)$ and $(G',x')$
are isomorphic.

\begin{Definition}
\label{defBS-limit} Given a sequence of random, finite graphs
$G_n$, let $k_0^{(n)}$ be a uniformly random vertex from $G_n$.
Following \cite{BS01}, we say that an infinite random, rooted
graph $(G,x)$ is the weak local limit of $G_n$ if for all
finite rooted graphs $(H,y)$ and all finite $r$, the
probability that $B_{r}(G_n,k_0^{(n)})$ is isomorphic to
$(H,y)$ converges to the probability that $B_{r}(G,x)$ is
isomorphic to $(H,y)$.
\end{Definition}

The main result of the paper is the following theorem.

\begin{Theorem}\label{thmmain}
The weak local limit of the all three variations of the
preferential attachment model is the P\'olya-point graph.
\end{Theorem}

Recently, and independently of our work, Rudas et al.
\cite{rudas}, studied the random tree resulting from the
preferential attachment model when $m = 1$. They derived the
asymptotic distribution of the subtree under a randomly
selected vertex which implies the Benjamini--Schramm limit.
Note that when $m=1$, there is no distinction between the
independent, conditional and sequential models.

As alluded to before, the points $x_{\bar a}$ of the P\'olya-point
process represent the random variables $S_{k-1}$ of the
vertices in $G_n$, which in turn behave like $(k/n)^\chi$ as
$n\to\infty$. The variable $y_{\bar a}=x_{\bar a}^{1/\chi}$
thus represents the birth-time of the corresponding vertex in
$G_n$. This is made precise in the following corollary to the
proof of Theorem \ref{thmmain}. As the theorem, the corollary
holds for all three versions of the Preferential Attachment
model.

\begin{Corollary}
\label{corlimit} Given $r<\infty$ and $\varepsilon>0$ there exists a
$n_0<\infty$ such that for $n\geq n_0$, there exists a coupling $\mu$
between a sample $T$ of the P\'olya-point and an ensemble $\{G_n,
v_0\}$ where $G_n$ has the distribution of the preferential attachment
graph of size $n$, and $v_0$ is a uniformly chosen vertex of $G_n$,
satisfying: with $\mu$ probability at least $1-\varepsilon$, there
exists an isomorphism $\bar a\mapsto k_{\bar a}$ from the ball of
radius $r$ about $0$ in $(T,0)$ into the ball of radius $r$ about $v_0$
in $G_n$, with the property that
\[
\biggl\llvert y_{\bar a}-\frac{k_{\bar a}}n\biggr\rrvert\leq\varepsilon
\]
for all $\bar a$ with distance at most $r$ from the root in
$(T,0)$. Here $y_{\bar a}$ is defined as $y_{\bar a}=x_{\bar
a}^{1/\chi}$.
\end{Corollary}

The numerator ${x_{\bar{a}}^{\psi}}=y_{\bar a}^{1-\chi}$ in
(\ref{poisson-intensity}) thus expresses the fact that in the
preferential attachment process, earlier vertices are more
likely to attract many neighbors than later vertices.\vadjust{\goodbreak}
%

\subsection{Subgraph frequencies}
\label{secsubgraphfrequency} A natural question concerning a
sequence of growing graphs $(G_n)$ is the question of how often
a small graph $F$ is contained in $G_n$ as a subgraph. This
question can be formalized in several ways, for example, by
considering the number of homomorphisms from $F$ into $G_n$, or
the number of injective homomorphism, or the number of times
$F$ is contained in $G_n$ as an \textit{induced} subgraph.

For graph sequences with bounded degrees, this leads to an
alternative notion of convergence, by defining sequence of
graphs to be convergent if the homomorphism density
$t(F,G_n)$---defined as the number of homomorphisms from $F$ into $G_n$
divided by the number of vertices in $G_n$---converges for all
finite connected graphs $F$ \cite{BCLSV-rev,sparse}. Indeed, for
sequences of graphs $G_n$ whose degree is bounded uniformly in
$n$, this notion can easily be shown to be equivalent to the
notion introduced by Benjamini and Schramm; moreover, the
corresponding notions involving the number of injective
homomorphisms, or the number of induced subgraphs, are
equivalent as well; see \cite{BCLSV-rev}, Section 2.2 for
formulas expressing these various numbers in terms of each
other.

But for graphs with growing maximal degree, this equivalence
does not hold in general. Indeed, consider a sequence of
graphs with uniformly bounded degrees, augmented by a vertex of
degree $n^{1/2}$. Such a vertex does not change the notion of
convergence introduced by Benjamini and Schramm; however, the number
of homomorphisms from a star with $3$ legs into this graph
sequence grows like $n^{3/2}$, implying that the homomorphism
density diverges.

To overcome this difficulty, we will consider maps $\Phi$ from
$V(F)$, the vertex set of $F$, into $V(G_n)$, the vertex set of
$G_n$ which in addition to being homormorphisms also preserve
degrees. More explicitly, given a graph $F$ and a map $
\mathbf{n}\dvtx V(F)\to\{0,1,2,\ldots\}$, we define $\inj(F,\mathbf{n};G_n)$ as
the number of injective maps $\Phi\dvtx V(F)\to V(G_n)$ such that:
\begin{longlist}[(2)]
\item[(1)] If $ij\in E(F)$, then $\Phi(i)\Phi(j)\in E(G_n)$;
\item[(2)] $d_{\Phi(i)}(G_n)=d_i(F)+n(i)$ for all $i\in V(F)$,
\end{longlist}
where $E(F)$ denotes the set of edges in $F$, and $d_i(F)$
denotes the degree of the vertex $i$ in $F$.

The following lemma is due to Laci Lovasz.

\begin{Lemma}
\label{lemsub-G-conv}
Let $D<\infty$, and let $G_n$ be a sequence of graphs
that
converges in the sense of Benjamini and Schramm. Then the
limit
%
\[
\hat t(F,\mathbf{n})=\lim_{n\to\infty}\frac1{|V(G_n)|}
\inj(F,\mathbf{n};G_n)
\]
%
exists for all finite connected graphs $F$ and all maps $
\mathbf{n}\dvtx V(F)\to\{0,1,2,\ldots\}$.
\end{Lemma}

As stated, the lemma refers to sequences of deterministic graphs. For sequences
of random graphs, its proof\vspace*{2pt} gives convergence of the expected number of
the subgraph frequencies $\frac1{|V(G_n)|}
\inj(F,\mathbf{n};G_n)$. To prove\vspace*{2pt} convergence in probability for these
frequencies, a little more work is needed.
For the case of preferential attachment graphs, we do this in
Section \ref{secfinite-ball}, together with an
explicit calculation of the actual values of these numbers.

\begin{Remark}
When $G_n$ has multiple edges, the definition of $
\inj(F,\mathbf{n};G_n)$ has to be modified. There are a priory
several possible definitions; motivated by the notions
introduced in \cite{BCLSV-rev} we chose the definition
\[
\inj(F,\mathbf{n};G_n)= \sum_{\Phi}\prod
_{ij\in E(F)}m_{\Phi(i)\Phi(j)}(G_n)^{m_{ij}(H)},
\]
where the sum goes over injective maps $\Phi\dvtx V(H)\to V(G_n)$
obeying condition (2) above with $d_i(H)$ and $d_{\Phi(i)}(G_n)$ denoting
degrees
counting multiplicities, and where $m_{ij}(H)$ is the multiplicity of
the edge $ij$ in $H$ [and similarly for $m_{\Phi(i)\Phi(j)}(G_n)$].
With this definition, the above lemma holds for graphs with
multiple edges as well.
\end{Remark}

\section{Proof of weak distributional convergence for the sequential
model}\label{secpf}
In this section we prove that the sequential model converges to
the P\'olya-point tree.

\subsection{P\'olya urn representation of the sequential model}
\label{secpolya}

In the early twentieth century, P\'olya proposed and analyzed the
following model known as the P\'olya urn model; see
\cite{durrett}. The model is described as follows. We have a
number of urns, each holding a number of balls, and at each
step, a new ball is added to one of the urns. The probability
that the ball is added to urn $i$ is proportional to $N_i + u$
where $N_i$ is the number of balls in the $i$th urn and $u$ is
a predetermined parameter of the model.

P\'olya showed that this model is equivalent to another process
as follows. For every $i$, choose at random a parameter (which
we call ``strength'' or ``attractiveness'')~$p_i$, and at each
step, \emph{independently} of our decision in previous steps,
put the new ball in urn $i$ with probability $p_i$. P\'olya
specified the distribution (as a function of $u$ and the
initial number of balls in each urn) for which this mimics the
urn model. A particularly nice example is the case of two urns,
each starting with one ball and $u=0$. Then $p_1$ is a uniform
$[0,1]$ variable, and $p_2=1-p_1$. P\'olya showed that for
general values of $u$ and $\{N_i(0)\}$, the values of $\{p_i\}$
are determined by the $\beta$-distribution with appropriate
parameters.

It is not hard to see that there is a close connection between
the preferential attachment model of Barab\'asi and Albert
and the P\'olya urn model in the following sense: every new
connection that a vertex gains can be represented by a new ball
added in the urn corresponding to that vertex.\vadjust{\goodbreak}


To derive this representation, let us consider first a two-urn
model, with the number of balls in one urn representing the
degree of a particular vertex $k$, and the number of balls in
the other representing the sum of the degrees of the vertices
$1,\ldots, k-1$. We will start this process at the point when
$n=k$ and $k$ has connected to precisely $m$ vertices in
$\{1,\ldots, k-1\}$. Note that at this point, the urn
representing the degree of $k$ has $m$ balls, while the other
one has $(2k-3)m$ balls.

Consider a time in the evolution of the preferential attachment
model when we have $n-1\geq k$ old vertices, and $i-1$ edges
between the new vertex $n$ and $\{1,\ldots, k-1\}$. Assume that
at this point the degree of $k$ is $d_k$, and the sum of the
degrees of $1,\ldots,k-1$ is $d_{<k}$. At this point, the
probability that the $i$th edge from $n$ to
$\{1,\ldots, n-1\}$ is attached to $k$ is
%
%
\begin{eqnarray}
\label{eqPntok}
&&\tilde\alpha\frac{1}{n-1} + (1-\tilde\alpha)
\frac{d_k}{2m(n-2) + (i-1)} \nonumber\\[-8pt]\\[-8pt]
&&\qquad=\frac{2m\alpha+ (1-\alpha
)d_k}{2m(n-2)+2m\alpha+(1-\alpha)(i-1)},\nonumber
\end{eqnarray}
while the probability that it is attached to one of the nodes
$1,\ldots, k-1$ is
%
%
\begin{eqnarray}
\label{eqPnto<k}
&&
\tilde\alpha\frac{k-1}{n-1} + (1-\tilde\alpha)
\frac
{d_{<k}}{2m(n-2) + (i-1)} \nonumber\\[-8pt]\\[-8pt]
&&\qquad=\frac{2m\alpha+ (1-\alpha
)d_{<k}}{2m(n-2)+2m\alpha+(1-\alpha)(i-1)}.\nonumber
\end{eqnarray}
Thus, conditioned on connecting to $\{1,\ldots,k\}$, the
probability that the $i$th edge from $n$ to
$\{1,\ldots, n-1\}$ is attached to $k$ is
\[
\frac1Z (2mu+d_k ),
\]
while the conditional probability that it is attached to one of
the nodes $1,\ldots, k-1$ is
\[
\frac1Z \bigl(2mu(k-1)+d_{<k} \bigr),
\]
where
$Z$ is an appropriate normalization constant. Note that the
constant $\tilde\alpha$ in (\ref{tilde-alpha}) was chosen in
such a way that the factor $u$ appearing in these expressions
does not depend on $i$, which is crucial to guaranty
exchangeability.

Taking into account that the two urns start with $m$ and
$(2k-3)m$ balls, respectively, we see that the evolution of the
two bins is a P\'olya urn with strengths $\psi_k$ and
$1-\psi_k$, where $\psi_k\sim
B_k=\beta(m+2mu,(2k-3)m+\break 2mu(k-1) )$.

\begin{pf*}{Proof of Theorem \ref{thm1}}
Using the two urn process as an inductive input, we can now
easily construct the P\'olya graph defined in
Theorem \ref{thm1}. Indeed, let $X_t\in\{1,2,\ldots,\lceil
\frac tm\rceil\}$ be the vertex receiving the $t$th
edge in the sequential model (the other endpoint of this edge
being the vertex $\lceil\frac tm\rceil+1$). For $t\leq m$,
$X_t$ is deterministic (and equal to $1$), but starting at
$t=m+1$, we have a two-urn model, starting with $m$ balls in
each urn. As shown above, the two urns can be described as
P\'olya-urns with strengths $1-\psi_2$ and $\psi_2$. Once
$t>2m$, $X_t$ can take three values, but conditioned on
$X_t\leq2$, the process continues to be a two-urn model with
strengths $1-\psi_2$ and $\psi_2$. To determine the
probability of the event that $X_t\leq2$, we now use the above
two-urn model with $k=3$, which gives that the probability of
the event $X_t\leq2$ is $1-\psi_3$, at least as long as $t\leq
3m$. Combining these two-urn models, we get a three-urn model
with strengths $(1-\psi_2)(1-\psi_3)$, $\psi_2(1-\psi_3)$ and
$\psi_3$. Again, this model remains valid for $t>3m$, as long
as we condition on $X_t\leq3$.

Continuing inductively, we see that the sequence $X_t$ evolves
in stages:
\begin{itemize}
\item For $t=1,\ldots,m$, the variable $X_t$ is
deterministic: $X_t=1$.
\item For $t=m+1,\ldots, 2m$, the distribution of
$X_t\in\{1,2\}$ is described by a two-urn model with
strengths $1-\psi_2$ and $\psi_2$, where $\psi_2\sim
B_2$.
\item In general, for $t=m(k-1)+1,\ldots, km$, the
distribution of $X_t\in\{1,\ldots,k\}$ is described by a
$k$-urn model with strengths
%
%
\begin{equation}
\label{phi-jk} \phi_j^{(k)}=\psi_j\prod
_{i=j+1}^k (1-\psi_i),\qquad j=1,\ldots, k.
\end{equation}
Here $\psi_k\sim B_k$ is chosen at the beginning of the
$k$th stage, independently of the previously
chosen strengths $\psi_1,\ldots,\psi_{k-1}$ (for
convenience, we set $\psi_1=1$).
\end{itemize}

Note that\vspace*{1pt} the random variables $\phi_j^{(k)}$ can be expressed
in terms of the random variables introduced in
Theorem \ref{thm1} as follows: by induction on $k$, it is
easy to show that
%
%
\begin{equation}
\label{S-k-prod} S_k=\prod_{j=k+1}^n(1-
\psi_k).
\end{equation}
This implies that
\[
\phi_j^{(k)}= \frac{\psi_j}{S_k},
\]
which relates the strengths $\phi_j^{(k)}$ to the random
variables defined in Theorem \ref{thm1}, and shows that the
process derived above is indeed the process given in the
theorem.
\end{pf*}

In order to apply Theorem \ref{thm1}, we will use two
technical lemmas, whose proofs will be deferred to a later
section. The first lemma states a law of large numbers for the
random variables $S_k$.

\begin{Lemma}\label{lemSk}
For every $\varepsilon$ there exist $K<\infty$ such that for
$n\geq K$, we have that with probability at least $1-\varepsilon$,
\[
\max_{k\in\{1,\ldots,n\}}\biggl\llvert S_k- \biggl(\frac kn
\biggr)^\chi\biggr\rrvert\leq\varepsilon
\]
and
\[
\max_{k\in\{K,\ldots,n\}}\biggl\llvert S_k- \biggl(\frac kn
\biggr)^\chi\biggr\rrvert\leq\varepsilon\biggl(\frac kn
\biggr)^\chi.
\]
\end{Lemma}

The second lemma concerns a coupling of the sequence
$\{\psi_k\}_{k\geq1}$ and an i.i.d. sequence of
$\Gamma$-random variables $\{\chi_k\}_{k\geq1}$, where
$\chi_k\sim\Gamma(m+2mu,1)$. To describe the coupling, we
define a sequence of functions $f_k\dvtx[0,\infty)\to[0,1)$ by
%
%
\begin{equation}
\label{coup1} \prob\bigl(\psi_k\leq f_k(x)\bigr)=\prob(
\chi_k\leq x).
\end{equation}
Then $f_k(\chi_k)$ has the same distribution as $\psi_k$,
implying that $ (\{\chi_k\}_{k\geq1},\break \{f_k(\chi_k)\}_{k\geq
1}) $ defines a coupling between $\{\chi_k\}_{k\geq1}$ and
$\{\psi_k\}_{k\geq1}$.

\begin{Lemma}\label{lemfk}
Let $f_k$ be as in (\ref{coup1}), and let $\{\chi_k\}_{k\geq
1}$ i.i.d. random variables with distribution
$\Gamma(m+2mu,1)$. Given $\varepsilon>0$ there exist a $K<\infty$
so that the following holds:

\begin{longlist}
\item
With probability at least $1-\varepsilon$,
%
%
\begin{equation}
\label{chi-k-bd} \chi_k\leq\log^2 k \qquad\mbox{for all } k\geq
K;
\end{equation}

\item For $k\geq K$ and $x\leq\log^2 k$,
%
%
\begin{equation}
\label{coupling-bd} \frac{1-\varepsilon}{ 2mk(1+u)} x \leq f_k(x) \leq
\frac{1+\varepsilon}{ 2mk(1+u)} x.
\end{equation}
\end{longlist}
\end{Lemma}
We defer the proof of Lemmas \ref{lemSk} and \ref{lemfk} to
Section \ref{secpolest}.

\subsection{The exploration tree of $G_n$}
\label{secexpltree}

Let $K_r=K_r(G_n,k_0)$ be the set of vertices in $G_n$ which
have distance at most $r$ from the random root $k_0$, and let
$\hat B_r(G_n,k_0)$ be the graph on $K_r$ that contains all
edges in $G_n$ for which at least 
one endpoint has distance
$\leq r$ from $k_0$. When proving that the preferential model
converges to the P\'olya-point graph, we will use the notion of
convergence given in Definition~\ref{defBS-limit}, but instead
of the standard ball of radius $r$, we will use the
modified ball $\hat B_r(G_n,k_0)$. (It is obvious that this
definition is equivalent.)



We will prove our results by induction on $r$, using the
exploration procedure outlined in Section \ref{secexplore} in
the inductive step. To this end, it will be convenient to endow
the rooted graph $(G_n,k_0)$ with a structure which is similar
to the one defined for the P\'olya-point graph.\vadjust{\goodbreak} More precisely, we
will inductively define a rooted tree $(T^{(n)}_r,0)$ on
sequence of integers $\bar{a}=(0,a_1,a_2,\ldots,a_l)$, and a
homomorphism
\[
{\mathbf k}^{(r)}\dvtx T^{(n)}_r\to\hat
B_r(G_n,k_0)\dvtx \bar a\mapsto
k_{\bar a}
\]
as follows.

We start our inductive definition by mapping $0$ into a vertex
$k_0$ chosen uniformly at random from the vertex set
$\{1,\ldots,n\}$ of $G_n$. Given a vertex $\bar
a=(0,a_1,a_2,\ldots,a_l)\in T^{(n)}_r$ and its image $k_{\bar
a}$ in $G_n$, let $d_{\bar a}$ be the degree of $k_{\bar a}$ in
$G_n$, and let $k_{\bar a_-},k_1,\ldots, k_{d_{\bar a}-1}$ be
the neighbors of $k_{\bar a}$ in $G_n$, where $\bar
a_-=(0,a_1,a_2,\ldots,a_{l-1})$. Recalling that edges were
created one by one during the sequential preferential
attachment process, we order $k_1,\ldots, k_{d_{\bar a}-1}$ in
such a way that for all $i=1,\ldots, d_{\bar a}-2$, the edge
$(k_{\bar a},k_i)$ was born before the edge $(k_{\bar
a},k_{i+1})$. We then define the children of $\bar a$ to be
the points $(\bar a,1),\ldots,(\bar a,d_{\bar a}-1)$. This
defines~$T_{r+1}^{(n)}$. The map ${\mathbf k}^{(r+1)}$ is the
extension of $\mathbf k^{(r)}$ which maps $(\bar
a,1),\ldots,(\bar a,d_{\bar a}-1)$ to the vertices $k_1,\ldots,
k_{d_{\bar a}-1}$, respectively. We call a vertex $(\bar a,i)$
early or of type { $L$} if $k_{i}<k_{\bar a_-}$ and late or of
type { $R$} otherwise. Note that the root and vertices of type
{ $L$} have $m$ children of type { $L$}, while vertices of type { $R$}
have $m-1$ children of type { $L$}.

To make the dependence on $G_n$ explicit, we often use the
notation $T_r(G_n)$ for the tree $T_r^{(n)}$, and the notation
$\mathbf k^{(r)}(G_n)$ for the map $\mathbf k^{(r)}$. Note that
$\mathbf k^{(r)}$ does not, in general, give a graph
isomorphism between $T_r^{(n)}$ and $\hat B_r(G_n,k_0)$. But if
the map is injective when restricted to $T_r^{(n)}$, it is a
graph isomorphism. To prove Theorem \ref{thmmain}, it is
therefore enough to show that for all $r$, the map $\mathbf
k^{(r)}$ is injective and the tree $T_r^{(n)}$ converges in
distribution to $T_r$, the ball of radius $r$ in the P\'olya-point
graph $(T,0)$.

\subsection{Regularity properties of the P\'olya-point process}

In order to prove Theorem \ref{thmmain}, we will use some
simple regularity properties of the P\'olya-point process.

Recall the definition of the P\'olya-point graph $(T,0)$ and the
P\'olya-point process $\{x_{\bar a}\}$ from Section \ref
{secpolyapointdef}, as well as the notation $\rho_{\bar a}(x)\,dx$ for
the intensity defined in (\ref{poisson-intensity}). As usual,
we define the height of a vertex $\bar
a=(0,a_1,a_2,\ldots,a_l)$ in $T$ as its distance $l$ from the
root. We denote the rooted subtree of height $r$ in $(T,0)$ by
$(T_r,0)$.
\begin{Lemma}
\label{lemreg1} Fix $0\leq r<\infty$ and $\varepsilon>0$. Then
there are constants $\delta>0$, \mbox{$C<\infty$}, $K<\infty$ and
$N<\infty$ such that with probability at least $1-\varepsilon$, we
have that:
\begin{itemize}
\item$x_{\bar a}\geq\delta$ for all vertices $\bar a$ in
$T_r$;
\item$\gamma_{\bar a}\leq C$;
\item$\rho_{\bar a}(\cdot)\leq K$;
\item$|T_r|\leq N$.\vadjust{\goodbreak}
\end{itemize}
\end{Lemma}

\begin{pf}
The proof of the lemma is easily obtained by induction on $r$.
We leave it to the reader.
%
%
%
\end{pf}

\begin{Corollary}
\label{correg2} For all $\varepsilon>0$ and all $r<\infty$ there
is a constant $\delta>0$ such that with probability at least
$1-\varepsilon$, we have
\[
\mathop{\min_{\bar a,\bar b\in T_r}}_{\bar a\neq\bar b} |x_{\bar
b}-x_{\bar
a}|\geq\delta.
\]
\end{Corollary}
\begin{pf}
This is an immediate consequence of the continuous nature of
the random variables $x_{\bar a}$ and the statements of
Lemma \ref{lemreg1}.
\end{pf}

\subsection{The neighborhood of radius one}
\label{sec1-Neighborhood}

Before proving our main theorem, Theorem \ref{thmmain}, for
the sequential model, we establish the following lemma, which
will serve as the base in an inductive proof of our main
theorem.

\begin{Lemma}
\label{thm2} Let $G_n$ be the sequential preferential
attachment graph, let $k_0$ be chosen uniformly at random in
$\{1,\ldots,n\}$ and let $k_{0,1},\ldots,k_{0,m+q_0}$ be the
neighbors of $k_0$, ordered as in Section \ref{secexpltree}
by the birth times of the edges $\{k_0,k_{0,i}\}$. Then
$(G_n,k_0)$ and the P\'olya-point process $\{x_{\bar a}\}$ can be
coupled in such a way that for all $\varepsilon>0$ there are
constants $C,N<\infty$, $\delta>0$ and $n_0<\infty$ such that
for $n\geq n_0$, with probability at least $1-\varepsilon$, we
have that:
\begin{longlist}[(iii)]
\item[(i)] $T_1 \cong T_1(G_n)$ and $|T_1(G_n)|\leq N$;
\item[(ii)] $\llvert x_{\bar a}- S_{k_{\bar a}-1}\rrvert\leq
\varepsilon$ for all $\bar a\in T_1$;
\item[(iii)] $k_0,k_{0,1},\ldots,k_{0,m+q_0}$ are pairwise distinct
and $k_{\bar a}\geq\delta n$ for all $\bar a\in T_1$;
\item[(iv)] $ \chi_{k_{\bar a}}=\gamma_{\bar a}\leq C$
for all $\bar a\in T_1. $
\end{longlist}
\end{Lemma}


\begin{pf}
(i)--(ii): We start by proving the first two statements. Choose
$y_0$ uniformly at random in $[0,1]$, let $x_0=y_0^{\chi}$ and
let $x_{0,1},\ldots,x_{0,m+q_0'}$ be the\vspace*{2pt} positions of the
children of $0$ in $(T,0)$. Define $k_0=\lceil ny_0\rceil$, so
that $k_0$ is distributed uniformly in $\{1,\ldots,n\}$, and for
$i=1,\ldots,m$, define $k_{0,i}$ by
\[
S_{k_{0,i}-1}\leq\frac{x_{0,i}}{x_0} 
S_{k_0-1} <
S_{k_{0,i}}.
\]
By\vspace*{2pt} Theorem \ref{thm1} and the observation that
$U_{k_0,1}=\frac{x_{0,1}}{x_0},\ldots,U_{k_0,m}=\frac{x_{0,m}}{x_0}$
are i.i.d. random variables
chosen uniformly at random from $[0,1]$, we have that indeed, with
large probability,
$k_{0,1},\ldots,k_{0,m}$ are close enough to the $x_{0,i}$'s.

Indeed,
given $\varepsilon>0$ choose $\delta$, $C$, $K$ and $N$ in such a
way that the statements of Lemma \ref{lemreg1} and
Corollary \ref{correg2} hold for $r=1$, and let
$\varepsilon'=\min\{\varepsilon,\delta/4\}$. By Lemma~\ref{lemSk}
there exists a constant $n_0<\infty$ such that for $n\geq n_0$,
we have that
%
%
\begin{equation}
\label{S-x-upbd1}\quad |S_{k_0-1}-x_0|\leq\varepsilon'
\quad\mbox{and}\quad |S_{k_{0,i}-1}-x_{0,i}| \leq\varepsilon'
\qquad\mbox{for all } i=1,\ldots,m
\end{equation}
with probability at least $1-2\varepsilon$.\vadjust{\goodbreak}

To understand the limiting distribution of the remaining
neighbors, $k_{0,m+1},\allowbreak\ldots, k_{0,m+q_0}$, of $k_0$, we observe
that conditioned on the random variables $\psi_1,\ldots,\psi_n$,
each vertex $k>k_0$ has $m$
independent
chances of being
connected to $k_0$, corresponding to the $m$ independent events
$\{X_{k,i}=k_0\}$, $i=1,\ldots,m$, where we used the shorthand
$X_{k,i}$ for the interval containing the endpoint of the $i$th edge
sent out
from $k$ (it is related to the random variables $X_t$
introduced in the proof of Theorem~\ref{thm1} via
$X_{k,i}=X_{(k-2)m+i}$). Let
%
%
\begin{equation}
\label{Pktok} P_{k\to k_0}=\phi_{k_0}
\frac1{S_{k-1}} =\frac{S_{k_0}}{S_{k-1}}\psi_{k_0}
\end{equation}
be the probability of the event $\{X_{k,i}=k_0\}$, and let
$N_{y_0}(y)=\break  \sum_{i=1}^m\sum_{k=k_0}^{\lceil ny\rceil}\mathbb
I(X_{k,i}=k_0)$ where $\mathbb I(A)$ is the indicator function
of the event $A$. We want to show that $N_{y_0}(\cdot)$
converges to a Poisson process on $[y_0,1]$.

By Lemma \ref{lemreg1}, we have that $k_0\geq nx_0\geq
n\delta$ with probability at least $1-\varepsilon$, which allows
us to apply Lemmas \ref{lemSk} and \ref{lemfk} to show that
for $n$ large enough, with probability at least $1-2\varepsilon$,
we have
\[
\hat P_{k\to k_0}(1-\varepsilon) \leq P_{k\to k_0} \leq(1+\varepsilon
)\hat
P_{k\to k_0}
\]
where
\[
\hat P_{k\to k_0} =\frac1{nm}
\frac{\chi_{k_0}}{2(1+u)}\frac n{k_0} \biggl(\frac{k_0}k
\biggr)^\chi.
\]
For $y>y_0$, let $\hat N_{y_0}(y)=
\sum_{i=1}^m\sum_{k=k_0}^{\lceil ny\rceil}\hat Y_{k\to
k_0}^{(i)}$ where $\{\hat Y_{k\to k_0}^{(i)}\}$ are independent
random variables such that $\hat Y_{k\to k_0}^{(i)}=1$ with
probability $\hat P_{k\to k_0}$ and $\hat Y_{k\to k_0}^{(i)}=0$
with probability $1-\hat P_{k\to k_0}$. It follows from
standard results on convergence to Poisson processes (and the
fact that $\gamma_0$ has the same distribution as $\chi_{k_0}$)
that $\hat N_{y_0}(\cdot)$ converges weakly to a Poisson
process with density $\frac{\gamma_0}{2(u+1)y_0} (\frac
{y_0}y )^\chi$ on $[y_0,1]$. A change of variables from $y$
to $x=y^\chi$ now leads to the Poisson process with density
\[
\frac{\gamma_0}{ 2(1+u)\chi}\frac{x^{\psi-1}}{x_0^\psi} { =\gamma_0
\frac{\psi x^{\psi-1}}{x_0^\psi}}
\]
on $[x_0,1]$. Combined with a last application of
Lemma \ref{lemSk} to bound the difference between
$S_{k_{0,i}-1}$ and $(k_{0,i}/n)^\chi$, this proves that
$x_{0,m+1},\ldots,x_{m+q_0'}\in[x_0,1]$ and
$k_{0,{m+1}},\ldots,k_{0,m+q_0}$ can be coupled in such a way
that for $n$ large enough, with probability at least
$1-3\varepsilon$, we have that $q_0=q_0'\leq Q=N-m-1$,
$\chi_{k_0}=\gamma_0\leq C$ and
%
%
\begin{equation}
\label{S-x-upbd2} |x_{0,i}-S_{k_{0,i}-1}|\leq\varepsilon'
\qquad\mbox{for } i=m+1,\ldots,m+q_0.
\end{equation}
Since $\varepsilon>0$ was arbitrary, this completes the proof of
the first two statements of the lemma.

(iii) To prove the third statement, we use bounds
(\ref{S-x-upbd1})
and (\ref{S-x-upbd2}), and a final
application of Lemma \ref{lemSk}, to establish the existence
of two constants $\delta'>0$ and $n_0'<\infty$ such that for
$n\geq n_0'$, with
probability
at least $1-4\varepsilon$,
%
%
\begin{equation}
\label{k-lbd1} k_{\bar a} \geq\delta'n \qquad\mbox{for all } \bar a
\in T_1(G_n)
\end{equation}
and
\[
|k_{\bar a}-k_{\bar b}|\geq\delta'n \qquad\mbox{for all }
\bar a,\bar b\in T_1(G_n) \mbox{ with } \bar a\neq\bar b,
\]
implying in particular that $k_0,k_{0,1},\ldots,k_{0,m+q_0}$ are
pairwise distinct.

(iv) To prove the last statement, let us assume that
$\gamma_0\leq C$, and that $k_{0,1},\ldots,\allowbreak k_{0,m+q}$ are
pairwise distinct, with $k_{0,i}<k_0$ for $i\leq m$,
$k_{0,i}>k_0$ for $i> m$, $\min k_{0,i}\geq n\delta'$ and $q\leq
Q$. Let $A$ be the event that we have chosen $k_0$ as the
uniformly random vertex and that the neighbors of $k_0$ are the
vertices $k_{0,1},\ldots,k_{0,m+q}$. Let $\chi^{A,\gamma_0}$ be
the collection of random variables $\{\chi_k\}_{k\neq k_0}$
conditioned on $\chi_{k_0}=\gamma_0$ and $A$.
We will want to
show that $\chi^{A,\gamma_0}$ can be coupled to a collection of
independent random variables $\{\hat\chi_k\}_{k\neq k_0}$ such
that $\chi^{A,\gamma_{0}}=\{\hat\chi_k\}_{k\neq k_0}$ with
probability at least $1-\varepsilon$, and
%
%
\begin{equation}
\label{psi-distr} \hat\chi_k \sim\cases{
F_k', &\quad if $k\in\{k_{0,1},\ldots,k_{0,m}\}$,
\cr
F_k, &\quad otherwise.}
\end{equation}

Let $\rho(\cdot\mid A,\chi_{k_0})$ be the density of the {
(multi-dimesional)}
random variable $\chi^{A,\gamma_0}$, and let $\prob(\cdot)$ be
the joint distribution of $G_n$ and the random variables
$\chi_1,\ldots,\chi_n$. By Bayes's theorem,
%
%
\begin{equation}
\label{Bayes} \rho(\cdot\mid A,\chi_{k_0}=\gamma_0)=
\frac{\prob(A\mid\cdot,\chi_{k_0}{ = \gamma_0})} {
\prob(A\mid\chi_{k_0}{ = \gamma_0})}\rho_0(\cdot),
\end{equation}
where $\rho_0$ is the original density of the random variables
$\{\chi_k\}_{k\neq k_0}$ (we denote the corresponding
probability distribution and expectations by $P_0$ and $E_0$,
resp.).

We thus have to determine the probability of $A$ conditioned on
$\chi_1,\ldots,\chi_n$. With the help of Theorem \ref{thm1},
this probability is easily calculated, and is equal to
\begin{eqnarray*}
\prob\bigl(A\mid\{\chi_k\}\bigr) 
&=&{ m!} \prod
_{i=1}^m P_{k_0\to k_{0,i}} \prod
_{j=1}^{q} mP_{k_{0,m+j}\to k_0} (1-P_{k_{0,m+j}\to k_0})^{m-1}
\\
&&{} \times\prod_{k>k_0:k\notin\{k_{0,m+1},\ldots,k_{0,m+q}\}} (1-P_{k\to
k_0}
)^m
\\
&=&{ m!}\prod_{i=1}^m
P_{k_0\to k_{0,i}} \prod_{j=1}^{q}
\frac{mP_{k_{0,m+j}\to k_0}}{1-P_{k_{0,m+j}\to k_0}} \prod_{k>k_0}
(1-P_{k\to k_0}
)^m,
\end{eqnarray*}
where $P_{k\to k'}$ is the conditional probability defined in (\ref
{Pktok}). By
Lemma \ref{lemSk}, this implies that given any $\varepsilon'>0$,
we can find $n_0<\infty$ such that for $n\geq n_0$, we have
that with probability at least $1-\varepsilon'$ with respect to
$P_0$,
\begin{eqnarray*}
&&\bigl(1{-\varepsilon'}\bigr)\prob\bigl(A\mid\{\chi_k\}
\bigr)
\\
&&\quad\leq{ m!} \Biggl(\prod_{i=1}^m
\psi_{k_{0,i}} \biggl(\frac{k_{0,i}}{k_0} \biggr)^\chi\prod
_{j=m+1}^{m+q}m\psi_{k_0} \biggl(
\frac{k_{0}}{k_{0,j}} \biggr)^\chi\Biggr) \exp\biggl( -m
\psi_{k_0}\sum_{k>k_0} \biggl(
\frac{{k_0}}{{k}} \biggr)^\chi\biggr)
\\
&&\quad\leq\bigl(1+{\varepsilon'}\bigr)\prob\bigl(A\mid\{
\chi_k\}\bigr).
\end{eqnarray*}
To estimate
$\prob(A\mid\chi_{k_0})=E_0[\prob(A\mid\{\chi_k\})]$, we
combined this bound with the deterministic upper bound
\begin{eqnarray*}
\hspace*{-4pt}&&\prob\bigl(A\mid\{\chi_k\}\bigr)
\\
\hspace*{-4pt}&&\quad\leq{ m!} \prod_{i=1}^m
P_{k_0\to k_{0,i}} \prod_{j=1}^{q}mP_{k_{0,m+j}\to k_0}
\leq\frac1n (m\psi_{k_0})^q\prod
_{i=1}^m \psi_{k_{0,i}}
\\
\hspace*{-4pt}&&\quad\leq
C'{m!} \Biggl(\prod
_{i=1}^m \psi_{k_{0,i}} \biggl(
\frac{k_{0,i}}{k_0} \biggr)^\chi\prod_{j=m+1}^{m+q}m
\psi_{k_0} \biggl(\frac{k_{0}}{k_{0,j}} \biggr)^\chi\Biggr)
\exp\biggl( -m\psi_{k_0}\sum_{k>k_0} \biggl(
\frac{{k_0}}{{k}} \biggr)^\chi\biggr),
\end{eqnarray*}
where $C'=(\delta')^{-(m+Q)}\sup_{n\geq1}e^{mn f_{\delta'
n}(C)}$.\vspace*{2pt}

These bounds imply that given any $\varepsilon'>0$, we can find an
$n_0<\infty$ such that for $n\geq n_0$, with probability at
least $1-\varepsilon'/2$ with respect to $P_0$, we have
\[
\sqrt{1-\varepsilon'}\prod_{i=1}^m
\frac{\psi_{k_{0,i}}}{E_0(\psi_{k_{0,i}})} \leq\frac{\prob(A\mid\{\chi
_k\})} {
\prob(A\mid\chi_{k_0})} \leq\sqrt{1+\varepsilon'}
\prod_{i=1}^m \frac{\psi_{k_{0,i}}}{E_0(\psi
_{k_{0,i}})}.
\]
With the help of Lemma \ref{lemfk}, this shows that for $n$
large enough, with probability at least $1-\varepsilon'$, we have
\[
\bigl(1-\varepsilon'\bigr)\prod_{i=1}^m
\frac{\chi_{k_{0,i}}}{E_0(\chi_{k_{0,i}})} \leq\frac{\prob(A\mid\{\chi
_k\})} {
\prob(A\mid\chi_{k_0})} \leq\bigl(1+\varepsilon'
\bigr)\prod_{i=1}^m \frac{\chi_{k_{0,i}}}{E_0(\chi_{k_{0,i}})}.
\]

Recalling (\ref{Bayes}) and the definition of the random
variables $\{\hat\chi_k\}_{k\neq k_0}$, we therefore have shown
that with probability at least $1-\varepsilon'$ with respect to~$P_0$,
%
%
\begin{equation}
\label{rho-rho1} \bigl(1-\varepsilon'\bigr)\hat\rho(\cdot) \leq\rho(
\cdot\mid A,\chi_{k_0}=\gamma_0) \leq\bigl(1+
\varepsilon'\bigr)\hat\rho(\cdot),
\end{equation}
where $\hat\rho$ is the density of the random variables
$\{\hat\chi_k\}_{k\neq k_0}$. (We denote the corresponding
product measure by $\hat P$.)

To continue, we need to transform statements which happen with
high probability with respect to $P_0$ into statements which
happen with high probability with respect to $\hat P$. To this
end, we consider the general case of two probability measures
$\mu$ and $\nu$\vadjust{\goodbreak} such that $\nu$ is absolutely continuous with
respect to $\mu$, $\nu=f\mu$ for some nonnegative function
$f\in L_2(\mu)$. Let $\Omega_0$ be an event which happens with
probability $1-\varepsilon'$ with respect to $\mu$. Then
%
%
\begin{equation}
\label{abs-cont} \nu\bigl(\Omega_0^c\bigr)=\int
f1_{\Omega_0^c}\leq\sqrt{E_\mu\bigl(f^2\bigr)
\mu\bigl(\Omega_0^c\bigr)} =\sqrt{
\varepsilon' E_\mu\bigl(f^2\bigr)},
\end{equation}
implying that $\Omega_0$ happens with probability at least
$1-\sqrt{\varepsilon' E_\mu(f^2)}$ with respect to~$\nu$.

Applying this bound to the probability measures $P_0$ and $\hat
P$, we see that bound (\ref{rho-rho1}) holds with
probability at least $1-\sqrt{2\varepsilon'}$ with respect to
$\hat P$, provided $n$ (and hence $k_{0,1},\ldots,k_{0,m}$)
is large enough. Using this fact, one then easily shows that
\[
\bigl\|\hat\rho-\rho(\cdot\mid A,\chi_{k_0}=\gamma_0)
\bigr\|_1 \leq2\varepsilon'+2\sqrt{2\varepsilon'}.
\]
Choosing $\varepsilon'$ sufficiently small
($\varepsilon'=\varepsilon^2/32$ is small enough), we see that the
right-hand side can be bounded by $\varepsilon$, which proves that
$\chi^{A,\gamma_0}$ and $\{\hat\chi_k\}_{k\neq k_0}$ can be
coupled in such a way that they are equal with probability at
least $1-\varepsilon$, as required.
\end{pf}

\subsection{Proof of convergence for the sequential model}

In this section we show that the sequential model converges to
the P\'olya-point graph. Indeed, we prove slightly more, namely the
following proposition:

\begin{Proposition}\label{propmain}
Given $\varepsilon>0$ and $r<\infty$, there are constants
$C,N<\infty$, $\delta>0$ and $n_0<\infty$ such that for $n\geq
n_0$, the rooted sequential attachment graph $(G_n,k_0)$ and
the P\'olya-point process $\{x_{\bar a}\}$ can be coupled in such a
way that with probability at least $1-\varepsilon$, the following
holds:
\begin{longlist}[(4)]
\item[(1)] \label{l1} $T_r(G_n)\cong T_r$ and $|T_r(G_n)|\leq N$;
\item[(2)]\label{l2} $|x_{\bar a}-S_{k_{\bar a}-1}|\leq
\varepsilon$ for all $\bar a\in T_r$;
\item[(3)]\label{l3} $\mathbf{k}^{(r)}(G_n)$ is injective, and
$k_{\bar a}\geq\delta n$ for all $\bar a\in T_r$;
\item[(4)]\label{l4} $\gamma_{\bar a}=\chi_{k_{\bar a}}\leq C$
for all $\bar a\in T_r$.
\end{longlist}
\end{Proposition}

\begin{pf}
For $r=1$, this follows from Lemmas \ref{thm2} and
\ref{lemreg1}.

Assume by induction that the lemma holds for $r<\infty$, and
fix $T_r$, $\mathbf k^{(r)}(G_n)$, $\{x_{\bar a}\}_{\bar a\in
T_r}$, $\{\gamma_{\bar a}\}_{\bar a\in T_r}$ and
$\{\chi_{k_{\bar a}}\}_{\bar a\in T_r}$ in such a way that
(1)--(4) hold (an event which has probability at
least $1-\varepsilon$ by our inductive assumption).

Consider a vertex $\bar a\in\partial T_{r}=T_r\setminus
T_{r-1}$. We want to explore the neighborhood of $k_{\bar a}$
in $G_n$. To this end, we note that for all $\bar b\in
T_{r-1}$, the neighborhood of $k_{\bar b}$ is already
determined by our conditioning on $\mathbf k^{(r)}(G_n)$,
implying in particular that none of the edges sent out from
$k_{\bar a}$ can hit a vertex $k\in K_{r-1}$, unless, of
course, $\bar a$ is of type { $R$}, and $k$ happens to be the
parent of $k_{\bar a}$---in which case the edge between $k$
and $k_{\bar a}$ is already present. To determine the children
of type { $L$} of the vertex $k_{\bar a}$, we therefore have to
condition on not hitting the set $K_{r-1}$.\vadjust{\goodbreak} But apart from
this, the process of determining the children of $k_{\bar a}$
is exactly the same as that of determining the children of the
root $k_0$. Since $|K_{r}|\leq N$, $k\geq\delta n$ for all
$k\in K_{r}$, and $\chi_{k}\leq C$ for all $k\in K_{r}$, we
have that $\sum_{k\in K_r} \phi_k\leq C'/n$ for some
$C'<\infty$, implying that conditioning on $k\notin
K_{r-1}\subset K_r$ has only a negligible influence on the
distribution of the children of $k_{\bar a}$. We may therefore
proceed as in the proof of Lemma \ref{thm2} to obtain a
coupling between a sequence of i.i.d. random variables $x_{\bar
a,i}$ distributed uniformly in $[0,x_{\bar a}]$ and the
children $k_{\bar a,i}$ of $k_{\bar a}$ that are of type ${ L}$.
As before, we obtain that for $n$ large enough, with
probability at least $1-\varepsilon$, we have $|S_{k_{\bar
a,i}-1}-x_{\bar a,i}|\leq\varepsilon$.

Repeating this process for all $k_{\bar a}\in\partial
K_{r}=K_r\setminus K_{r-1}$, we obtain a set of vertices
${ L}_{r+1}$ consisting of all children of type ${ L}$ with parents
in $\partial K_{r}$. It is easy to see that with probability
tending to one as $n\to\infty$, the set ${ L}_{r+1}$ has no
intersection with $K_r$, so we will assume this for the rest of
this proof.

Next we continue with the vertices of type $R$. Assume that we
have already determined all children of type $R$ for a certain
subset $U_r\subset\partial K_{r}$, and denote the set children
obtained so far by $R_{r+1}$. We decompose this set as
$R_{r+1}=\bigcup_{i=1}^m R_{r+1}^{(i)}$, where
$R_{r+1}^{(i)}=\{k\in R_{r+1}\dvtx X_{i,k}\in U_r\}$.

Consider a vertex $\bar a\in\partial K_r\setminus U_r$.
Conditioning on the graph explored so far is again not
difficult, and now amounts to two conditions:

\begin{longlist}[(2)]
\item[(1)] $X_{k,i}\neq k_{\bar a}$ if $k\in K_r\cup
R_{r+1}^{(i)}$,
since all the edges sent out from this set have already
been determined.

\item[(2)] For $k\notin K_r\cup R_{r+1}^{(i)}$, the probability
that $k_{\bar a}$ receives the $i$th edge
from $k$ is different from the probability given in
(\ref{Pktok}), since the random variables $X_{k,i}$
has been probed before: we know that $X_{k,i}\notin
K_{r-1}$ since otherwise $k$ had sent out an edge to a
vertex in $K_{r-1}$, which means that $k$ would have
been a child of type ${ R}$ in $K_r$. We also know that
$X_{k,i}\notin U_{r}$, since otherwise $k\in
{ R}_{r+1}^{(i)}$. Instead of (\ref{Pktok}), we therefore
have to use the modified probability
\[
P_{k\to k_{\bar a}}= \phi_{k_{\bar a}}\frac1{\tilde S_{k-1}},
\]
where
\[
\tilde S_{k-1} =\mathop{\sum_{k'>k_{\bar a}:}}_{k'\notin K_{r-1}\cup U_r}
\phi_{k'}.
\]
\end{longlist}
Since $\tilde S_{k-1}\leq S_{k-1}\leq\tilde S_{k-1}+C'/n$ by
our inductive assumption, we can again refer to
Lemma \ref{lemSk} to approximate $P_{k\to k_{\bar a}}$ by
\[
\hat P_{k\to k_{\bar a}} =\frac1{nm}\frac{\chi_{k_{\bar
a}}}{2(1+u)}\frac n{k_{\bar a}}
\biggl(\frac{k_{\bar a}}k \biggr)^\chi.
\]
From here on the proof of our inductive claim is completely
analog to the proof of Lemma \ref{thm2}. We leave it to the
reader to fill in the (straightforward but slightly tedious)
details.
\end{pf}

\subsection{Estimates for the P\'olya urn
representation}\label{secpolest}

In this section we complete the work started in Section \ref
{secpolya} by proving Lemmas \ref{lemSk} and
\ref{lemfk}.

\begin{pf*}{Proof of Lemma \ref{lemSk}}
Fix $\eps$, and recall that
\[
\chi=\frac{1+2u}{2+2u}\in\biggl[\frac12,1 \biggr).
\]
Writing $S_k$ as
\[
S_k=\prod_{j=k+1}^n(1-
\psi_j) =\exp\Biggl( \sum_{j=k+1}^n
\log(1-\psi_j) \Biggr),
\]
we use the fact that if $0<x<1$, then $x\leq-\log(1-x)\leq x
+x^2/(1-x)$ to
bound
\[
\Biggl\llvert E \Biggl[\sum_{j=k+1}^n
\log(1-\psi_j) \Biggr]+\sum_{j=k+1}^n
E[\psi_j]\Biggr\rrvert\leq\sum_{j=k+1}^n
E \biggl[\frac{\psi_j^2}{1-\psi_j} \biggr].
\]
On the other hand, by Kolmogorov's inequality and the fact that
\[
\operatorname{Var}\bigl(\log(1-\psi_j)\bigr)\leq E\bigl[\bigl(\log(1-
\psi_k)\bigr)^2\bigr]\leq E\bigl[\psi_j^2
(1-\psi_j)^{-2}\bigr],
\]
we have
\begin{eqnarray*}
&&
\prob\Biggl(\max_{K\leq k\leq n} \Biggl\llvert\sum
_{j=k+1}^n\log(1-\psi_j)-E \Biggl[\sum
_{j=k+1}^n\log(1-\psi_j)
\Biggr]\Biggr\rrvert\geq\varepsilon\Biggr)
\\
&&\qquad\leq\frac1{\varepsilon^2}\sum_{j=K+1}^n
E \biggl[\frac{\psi
_j^2}{(1-\psi_j)^2} \biggr].
\end{eqnarray*}

We will use that for any $\beta_{a,b}$ distributed random
variable $\psi$, we have
\[
E[\psi]=\frac a{a+b},\qquad E \biggl[\frac{\psi^2}{1-\psi} \biggr] =\frac
{a(a+1)}{(a+b)(b-1)}
\]
and
\[
E \biggl[\frac{\psi^2}{(1-\psi)^2} \biggr] =\frac
{a(a+1)}{(b-2)(b-1)}.
\]
Using these relations for $a=m+2mu$ and $b=(2j-3)m+2mu(j-1)$,
we get
%
%
\begin{eqnarray}
\label{eqexpepsi-0} E(\psi_j)&=& \frac{m+2mu}{(2j-2)m+2jmu}=
\frac{\chi}{j}+O \biggl(\frac1{j^2} \biggr),
\\
%
\label{beta-moments} E\bigl[\psi_j^2\bigr]&\leq& E
\biggl[\frac{\psi_j^2}{1-\psi_j} \biggr]=O \biggl(\frac1{j^2} \biggr)
\quad\mbox{and}\quad E \biggl[\frac{\psi_j^2}{(1-\psi_j)^2} \biggr]=O \biggl
(\frac1{j^2}
\biggr).
\end{eqnarray}
Putting these bounds together, and observing that
$\sum_{j=k+1}^n\frac1j=\log(n/k)+O(k^{-1})$, we get that
there exists a constant $K(\eps)$ not depending on $n$ such
that with probability at least $1-\varepsilon$, we have that
\[
\biggl(\frac{k}{n} \biggr)^\chi e^{-\varepsilon}<S_k<
\biggl(\frac{k}{n} \biggr)^\chi e^\varepsilon\qquad{\mbox{for
all } K(\eps)\leq k\leq n.}
\]
For $k<K(\eps)$, we bound $S_k\leq S_K$ to conclude that with
probability at least $1-\varepsilon$,
\[
\biggl\llvert S_k- \biggl(\frac kn \biggr)^\chi\biggr
\rrvert=O \biggl( \biggl(\frac Kn \biggr)^\chi\biggr).
\]
The lemma now follows.
\end{pf*}

\begin{pf*}{Proof of Lemma \ref{lemfk}}
(i) Let $a=m+2mu$, so that $\chi_k\sim\Gamma(a,1)$. Then
\[
\prob\bigl(\chi_k\geq\log^2k\bigr)\leq E
\bigl[e^{\chi_k/2}\bigr] e^{-
(\log^2 k)/2} =2^a k^{-(\log k)/2}.
\]
Since the right-hand side is sumable, this implies the first
statement of the lemma through the Borel--Cantelli lemma.

(ii) Let $b_k=(2k-3)m+2mu(k-1)-1$, and let
$\chi_k'=\chi_k/b_k$. Then $f_k$ can be defined by
\[
\prob\bigl(\psi_k\leq f_k(x)\bigr)=\prob\bigl(
\chi_k'\leq x/b_k\bigr).
\]
In order to prove the second statement of the lemma, it is
clearly enough to prove that for all sufficiently large $k$, we
have
\[
(1-\varepsilon)\frac x{b_k} \leq f_k(x) \leq\frac
x{b_k} \qquad\mbox{for } x\leq\log^2k,
\]
which in turn is equivalent to showing that
%
%
\begin{equation}
\label{proof-equ} \prob\bigl(\psi_k\leq(1-\varepsilon)x \bigr)\leq
\prob\bigl(\chi_k'\leq x \bigr)
\leq\prob(\psi_k\leq x ) \qquad\mbox{for } x\leq
\frac{\log^2k}{b_k}
\end{equation}
provided $k$ is large enough.

We start by proving that
\[
\Delta(x):=\prob(\psi_k\leq x)-\prob\bigl(\chi_k'
\leq x\bigr)\geq0.
\]
To this end, we rewrite
\[
\prob(\psi_k\leq x) =\frac1{Z_\beta}\int
_0^x y^{a-1}(1-y)^b\,dy
\]
and
\[
\prob\bigl(\chi_k'\leq\lambda\bigr) =
\frac1{Z_\gamma}\int_0^\lambda
y^{a-1}e^{-by}\,dy,
\]
where $a=m+2mu$, $b=b_k$ and $Z_\gamma=\int_0^\infty
y^{a-1}e^{-by}\,dy$ and $Z_\beta=\int_0^1y^{a-1}(1-y)^b\,dy$ are
the appropriate normalization factors. For $x\leq1$, we
express $\Delta(x)$ as
\[
\Delta(x)=\frac1{Z_\gamma}\int_0^x
dyy^{a-1}e^{-b y} \Biggl(e^\delta\exp\Biggl(-b\sum
_{k=2}^\infty\frac{y^k}k \Biggr)-1
\Biggr),
\]
where $e^\delta=Z_\gamma/Z_\beta$. Note that $\delta>0$ by the
fact that $(1-x)\leq e^{-x}$. It is also easy to see that
$\delta\to0$ as $k\to\infty$; indeed, we have
$\delta=O(b^{-1})=O(k^{-1})$.

Consider the derivative
\[
\frac{d\Delta(x)}{dx} =\frac{x^{a-1}e^{-b x}}{Z_\gamma} \Biggl(e^\delta
\exp\Biggl(-b
\sum_{k=2}^\infty\frac{x^k}k \Biggr)-1
\Biggr),
\]
and let $x_0$ be the unique root, that is, let $x_0\in(0,1)$ be
the solution of the equation
\[
\delta=b\sum_{k=2}^\infty
\frac{x_0^k}k.
\]
Then $\Delta(x)$ is monotone increasing for $0<x<x_0$ and
monotone decreasing for all $x>x_0$. Together with the
observation that $\Delta(x)>0$ for all sufficiently small $x$,
and $\Delta(x)\to0$ as $x\to\infty$, we conclude that
$\Delta(x)\geq0$ for $0\leq x<\infty$. This proves that
$\prob(\chi_k'\leq x)\leq\prob(\psi_k\leq x))$ for all $x\geq
0$.

To prove the lower bound in (\ref{proof-equ}), we will prove
that
\[
\tilde\Delta(x)=\prob\bigl(\chi_k'\leq x \bigr) -
\prob\bigl(\psi_k\leq(1-\varepsilon)x \bigr)\geq0 \qquad\mbox{if } x\leq\frac
\varepsilon4\leq\frac1{8}.
\]
We decompose the range of $x$ into two regions, depending on
whether $x\geq\frac{4a}{b\varepsilon}$ or $x\leq
\frac{4a}{\varepsilon b}$.

In the first region, we express $\tilde\Delta(x)$ as
\begin{eqnarray*}
\tilde\Delta(x) &=& \prob\bigl(\psi_k\geq(1-\varepsilon)x \bigr)-\prob
\bigl(\chi_k'\geq x \bigr)
\\
&=&\frac{e^\delta}{Z_\gamma}\int_{x(1-\varepsilon)}^1
dyy^{a-1}(1-y)^b -\frac1{Z_\gamma}\int
_x^\infty dyy^{a-1}e^{-b y}.
\end{eqnarray*}
We then bound
\begin{eqnarray*}
\int_{x}^\infty dy(2y)^{a-1}e^{-b y}
&\leq&\int_x^{2x} dyy^{a-1}e^{-b y}
\int_{2x}^\infty dyy^{a-1}e^{-b y}
\\
&\leq&\int_x^{2x} dyy^{a-1}e^{-b y}
+2^{a-1}e^{-bx}\int_{x}^\infty
dyy^{a-1}e^{-b y}
\end{eqnarray*}
proving that
%
%
\begin{eqnarray}
\label{Csbd1} \int_{x}^\infty dy(2y)^{a-1}e^{-b y}
&\leq&\bigl(1-2^{a-1}e^{-bx}\bigr)^{-1}\int
_x^{2x} dyy^{a-1}e^{-b y} \nonumber\\[-8pt]\\[-8pt]
&\leq&2
\int_x^{2x} dyy^{a-1}e^{-b y},\nonumber
\end{eqnarray}
where we have used $bx\geq a\log2$ in the last step.

On the other hand, using that $(1-y)^b\geq e^{-by(1+x)}$ if
$y\leq2x\leq1/2$, we have that
\begin{eqnarray*}
{e^\delta}\int_{x(1-\varepsilon)}^1
dyy^{a-1}(1-y)^b &\geq& \int_{x(1-\varepsilon)}^{2x(1-\varepsilon)}
dyy^{a-1}e^{-by(1+x)}
\\
&=& \int_{x}^{2x} dyy^{a-1}(1-
\varepsilon)^ae^{-by(1+x)(1-\varepsilon)}
\\
&\geq& (1-\varepsilon)^ae^{-2bx^2}e^{\varepsilon bx} \int
_x^{2x}dyy^{a-1}e^{-b y}
\\
&\geq& 2\int_x^{2x}dyy^{a-1}e^{-b y}.
\end{eqnarray*}
Combined with (\ref{Csbd1}), this proves that
$\tilde\Delta(x)\geq0$ if $\varepsilon b x\geq4a$.

For $\varepsilon b x\leq4a$, we bound
\begin{eqnarray*}
\tilde\Delta(x) &=& \frac1{Z_\gamma} \biggl(\int_0^x
dyy^{a-1}e^{-b y} -e^\delta\int_0^{x(1-\varepsilon)}dyy^{a-1}(1-y)^b
\biggr)
\\
&\geq& \frac1{Z_\gamma} \biggl(\int_0^x
dyy^{a-1}e^{-b y} -e^\delta\int_0^{x(1-\varepsilon)}dyy^{a-1}e^{-by}
\biggr)
\\
&=& \frac1{Z_\gamma} \biggl( \int_{(1-\varepsilon)x}^x
dyy^{a-1}e^{-b y} -\bigl(e^\delta-1\bigr)\int
_0^{x(1-\varepsilon)}dyy^{a-1}e^{-by}
\biggr)
\\
&\geq& \frac1{Z_\gamma} \bigl(\varepsilon x \bigl[(1-\varepsilon)x
\bigr]^{a-1}e^{-bx} -\bigl(e^\delta-1
\bigr)x^a \bigr)
\\
&\geq& \frac{x^a}{Z_\gamma} \bigl(\varepsilon2^{1-a}e^{-4a/\varepsilon} -
\bigl(e^\delta-1\bigr) \bigr).
\end{eqnarray*}
Since $\delta\to0$ as $b\to\infty$, we see that the right-hand
side becomes positive if $k\geq K$ for some $K<\infty$ that
depends on $a$ and $\varepsilon$ (it grows exponentially in
$a/\varepsilon$).
\end{pf*}

\section{Approximating coupling for the independent and the
conditional models}

In this section we prove that the sequential and the
independent model have the same\vadjust{\goodbreak} weak limit. To this end we
construct a coupling between the two models such with
probability tending to $1$, the balls around a randomly chosen
vertex in $\{1,\ldots,n\}$ are identical in both models. This
will imply that both models have the same weak local
limit.

We only give full details for the coupling between the
independent and the sequential model. The approximating
coupling between the conditional and the sequential model is
very similar, and the proof that it works is identical.

We construct the coupling inductively as follows: let
$V=1,2,\ldots$ be the vertices of the preferential attachment
graph. For $1\neq n\in V$ and $i=1,\ldots,m$ let $e^i_n<n$ and
$f^i_n<n$ be the $i$th vertex that $n$ is connected to in,
respectively, the sequential and the independent models. We
use the symbol $\mathbf e_n$ to denote the vector
$\{e^i_n\}_{1\leq i\leq m}$,
and the symbol $\mathbf
f_n$ to denote the vector $\{f^i_n\}_{i=1}^m$.

By construction, $e^i_2=f^i_2=1$ for all $i$. Once we know
$\mathbf e_l$ and $\mathbf f_l$ for every $l<n$, we determine
$\mathbf e_n$ and $\mathbf f_n$ as follows: let $D_1$ be the
distribution of $\mathbf e_n$, based on the sequential rule and
conditioned on $\{\mathbf e_l\}_{l<n}$, and let $D_2$ be the
distribution of $\mathbf f_n$ based on the independent rule
and conditioned on $\{\mathbf f_l\}_{l<n}$. Let $D$ be an (arbitrarily chosen)
coupling of $D_1$ and $D_2$ that minimizes the total variation
distance. Then we choose $\mathbf e_n$ and $\mathbf f_n$
according to $D$.

Our goal is to prove the following proposition:

\begin{Proposition}
\label{propcoupling} Let $(G_n)$ and $(G_n')$ be the sequence
of preferential attachment graphs in the sequential and the
conditional model, respectively, coupled as above. Let
$\varepsilon>0$ and let $r$ be an arbitrary positive integer.
Then there exists an integer $n_0$ such that for $n\geq n_0$,
with probability at least $1-\varepsilon$, a uniformly chosen random
vertex $k_0\in\{1,\ldots,n\}$ has the same $r$-neighborhood in
$G_n$ and $G_n'$.
\end{Proposition}

The proof of the proposition relies on following two lemmas, to
be proven
in 
Sections \ref{seccoupgood} and \ref{seclem-dk-proof},
respectively.

\begin{Lemma}\label{claimcoupgood}
Consider the coupling defined above, and fix $k\geq2$. For
$n>k$, let $A_n=A_n^{(k)}$ be the event that there exists an
$i\in\{1,\ldots,m\}$ such that $e^i_n=k\neq f^i_n$ or $e^i_n\neq
k=f^i_n$. Then
%
%
\begin{equation}
\label{eq1} \prob\Biggl(A_n\biggm|\bigcap_{h=k+1}^{n-1}A_h^c,d_{n-1}(k)
\Biggr) 
=O \biggl(\frac{d_{n-1}(k)}{n^2} \biggr).
\end{equation}
\end{Lemma}
Note that under the conditioning, $d_{n-1}(k)$ is the same in both models.

\begin{Lemma}\label{lemdk}
For the sequential preferential attachment model, for every $n$
and $k$ such that $n>k$, let $d_n(k)$ be the degree of
vertex $k$ when the graph contains $n$ vertices. Then
%
%
\begin{equation}
\label{eqdexp} E\bigl[d_n(k)\bigr]= m \biggl[1+\frac{\chi}{1-\chi}
\biggl( \biggl(\frac{n}{k} \biggr)^{1-\chi}-1 \biggr) \biggr] +O
\biggl(\frac{n^{1-\chi}}{k^{2-\chi}} \biggr),
\end{equation}
where the constant implicit in the $O$-symbol depends on $m$
and $u$.\vadjust{\goodbreak}
\end{Lemma}

\subsection{\texorpdfstring{Proof of Proposition \protect\ref{propcoupling}}
{Proof of Proposition 4.1}}

Fix $\varepsilon$ and $r$, let $B_r(k)$ and $B_r(k)'$ be the ball
of radius $r$ about $k$ in $G_n$ and $G_n'$, respectively, and
let $B$ be the set of vertices $k\in\{1,\ldots,n\}$ for which
$B_r(k)\neq B_r(k)'$. Then the probability that a uniformly
chosen vertex in $\{1,\ldots,n\}$ is in $B$ is just $1/n$ times
the expected size of $B$. We thus have to show that
\[
E\bigl[|B|\bigr]\leq\varepsilon n.
\]

In a preliminary step note that $B_r(k)=B_r(k)'$ unless there
exists a vertex $k'\in B_r(k)$ such that $e^i_{n'}=k'\neq
f^i_{n'}$ or $e^i_{n'}\neq k'=f^i_{n'}$ for some $i=1,\ldots, m$
and some $n'>k'$.

To prove this fact, let us consider the event
$A^{(k)}=\bigcup_{n>k} A_n^{(k)}$. It is easy to see that this
event is the event that at least one of the edges received by
$k$ is different in $(G_n)$ and $(G_n')$. Using this fact, one
easily shows that the ball of radius $1$ around a vertex $k$
must be identical in $G_n$ and $G_n'$ unless $A^{(k')}$ happens
for at least one vertex $k'$ in the $1$-neighborhood of $k$ in
$G_n$. By induction, this implies that $B_r(k)=B_r(k)'$ unless
there exists a vertex $k'\in B_r(k)$ such that the event
$A^{(k')}$ happens, which is what we claimed in the previous
paragraph.

Next we note that by Proposition \ref{propmain}, there exist
$\delta>0$ and $N<\infty$ such that with probability at least
$1-\varepsilon/2$, a random vertex $k\in\{1,\ldots,n\}$
obey the two following two conditions:
\begin{longlist}[(2)]
\item[(1)]\label{A1} the ball of radius $2r$ around $k$ in the
sequential graph $G_n$ contains no more than $N$
vertices;
\item[(2)]\label{A2} the oldest vertex (the vertex with the
smallest index) in this ball is no older than $\delta
n$.
\end{longlist}
If we denote the set of vertices satisfying these two
conditions by $W$, we thus have that
\[
E\bigl[|W|\bigr]\geq\biggl(1-\frac\varepsilon2\biggr)n.
\]
As a consequence, it will be enough to show that
\[
E\bigl[|W\cap B|\bigr]\leq\frac\varepsilon2 n.
\]
If $k\in W\cap B$, there must be a vertex $k'\in B_r(k)$ such
that the event $A^{(k)}$ happens. But $k'\in B_r(k)$ if and
only if $k\in B_r(k')$, and since $B_r(k')\subset B_{2r}(k)$,
we must further have that $|B_r(k')|\leq N$ and $k'\geq\delta
n$. As a consequence,
\begin{eqnarray*}
|W\cap B| &=& \sum_{k\in W}\I(k\in B) \leq\sum
_{k\in W}\sum_{k'\in B_r(k)}\I
\bigl(A^{(k')}\bigr)
\\
&=&\sum_{k'}\I\bigl(A^{(k')}\bigr)\sum
_{k\in B_r(k')}I(k\in W)
\\
&\leq& N\sum_{k'= \delta n}^{{n}}\I
\bigl(A^{(k')}\bigr),
\end{eqnarray*}
where we used the symbol $\I(A)$ to denote the indicator
function of the event~$A$.

Finally by Lemmas \ref{claimcoupgood} and \ref{lemdk},
\[
P\bigl(A^{(k)}\bigr) \leq O(1)\sum_{n>k}
\frac1{n^2} \biggl(\frac n k \biggr)^{1-\chi} 
= O
\biggl(\frac1k \biggr).
\]
As a consequence we can find a constant $C$ such that
\[
E\bigl[|W\cap B|\bigr] \leq N\sum_{k'= n\delta}^{{n}}\frac
C{k'} \leq CN/{\delta}.
\]
For $n$ large enough, the right-hand side is smaller than
$\frac\varepsilon2 n$, which is the bound
we had to establish.%
\subsection{\texorpdfstring{Proof of Lemma \protect\ref{claimcoupgood}}
{Proof of Lemma 4.2}}
\label{seccoupgood}

Let us the shorthand $d$ for the degree $d_{n-1}^k$.
In the
independent model the probability of having $r$ connections to
$k$ and $h=m-r$ connections to other vertices in
$\{1,\ldots,n-1\}$ is
\[
\pmatrix{m \cr r} p^r (1-p)^h \qquad\mbox{with } p=\frac
\alpha{n-1}+\frac{(1-\alpha)d}{2m(n-2)},
\]
while in the sequential model it is
\[
\pmatrix{m \cr r} \prod_{l=0}^{r-1}p_l
\prod_{l=r}^{m-1}(1-p_l)
\]
with
\[
p_l=p_l(r) =\cases{\displaystyle \frac{2m \alpha+
(1-\alpha)(d+l)}{2m(n-2)+2m\alpha+(1-\alpha)l}, &\quad if $l<r$,
\vspace*{2pt}\cr
\displaystyle \frac{2m \alpha+
(1-\alpha)(d+r)}{2m(n-2)+2m\alpha+(1-\alpha)l}, &\quad if $l\geq r$.}
\]
[Here we used exchangeability and (\ref{eqPntok}).]

As a consequence, the probability in (\ref{eq1}) is bounded by
a constant times
%
%
\begin{equation}
\label{eqrandj} \max_{r=0,\ldots,m} \Biggl\llvert\Biggl[ \prod
_{l=0}^{r-1}p_l\prod
_{l=r}^{m-1}(1-p_l) \Biggr] -
p^r (1-p)^h \Biggr\rrvert.
\end{equation}
Telescoping the difference, we bound (\ref{eqrandj}) by
\begin{eqnarray*}
&&\max_{r=0,\ldots,m} \Biggl( \sum_{l=0}^{r-1}p^l|p_l-p|
\prod_{l'=l+1}^{r-1}p_{l'} \\
&&\hspace*{17pt}\qquad{}+ \sum
_{l=r}^{m-1}p^r(1-p)^{l-r}\bigl|(1-p)-(1-p_l)\bigr|
\prod_{l'=l+1}^{m-1}(1-p_{l'})
\Biggr)
\\
&&\qquad\leq\max_{r=0,\ldots,m} \Biggl( {\tilde p}^{r-1}\sum
_{l=0}^{r-1}|p_l-p| +
p^r\sum_{l=r}^{m-1}|p-p_l|
\Biggr),
\end{eqnarray*}
where ${\tilde p}=\max\{p,p_1,\ldots,p_m\}=O(d/n)$. We now
distinguish three cases:

\begin{longlist}[(iii)]
\item[(i)] if $r\geq2$, we use the fact that $p-p_l=O(1/n)$ to get a
bound of order $ O({\tilde p}/n)=O(d/n^2)$ for both sums;

\item[(ii)] if $r=1$, we use the fact that the first sum is equal to
$|p_0-p|=O(1/n^2)$, while the second can be bounded by $
O({\tilde p}/n)=O(d/n^2)$ as before;

\item[(iii)]  if $r=0$, we use that fact that
\begin{eqnarray*}
p_l(0)&=&\frac{2m \alpha+ (1-\alpha)d}{2m(n-2)+2m\alpha+(1-\alpha)l}
\\
&=&\frac{2m \alpha}{2m(n-1)} \bigl(1+O\bigl(n^{-1}\bigr) \bigr) +
\frac{ (1-\alpha)d}{2m(n-2)} \bigl(1+O\bigl(n^{-1}\bigr) \bigr)
\\
&=&p+O\bigl(d/n^2\bigr)
\end{eqnarray*}
to show that for $r=0$, all terms in the sum
$\sum_{l=0}^{m-1}|p- p_l|$ are of order $O(d/n^2)$.
\end{longlist}

This completes the proof of the lemma.
\subsection{\texorpdfstring{Proof of Lemma \protect\ref{lemdk}}{Proof of Lemma 4.3}}
\label{seclem-dk-proof}

As before, we use $\phi_k^{(n)}$ for
\[
\phi_k^{(n)}=\psi_k\prod
_{i=k+1}^{n}(1-\psi_i).
\]
By construction,
%
%
\begin{equation}
\label{d-alt-expression} d_n(k)=m+\sum_{t=(k-1)m+1}^{(n-1)m}
\U_t,
\end{equation}
where\vspace*{1pt} the variables $\{\U_t\}$ are defined as follows: let
$\{\hat{U}_t\}_{t=1}^\infty$ be i.i.d. $U[0,1]$ variables,
independent of the $\phi_k$'s. Then
$
\U_t={\mathbf1}_{\hat{U}_t<\phi_k^{(\lceil t/m\rceil)}}.
$
Note that conditioned on
$\{\phi_k^{(j )}\}_{j\geq k}$, $\{\U_t\}$'s are independent,
each being Bernoulli $\phi_k^{(\lceil t/m\rceil)}$.

Let
$\FF$ be the $\sigma$-algebra generated by
$\{\psi_h\}_{h=1}^\infty$. Then
%
%
\begin{equation}
\label{eqhatd} E\bigl(d_n(k)\mid\FF\bigr) = m+m\sum
_{\ell= k}^{n-1} \phi_k^{(\ell)}.
\end{equation}
By (\ref{eqexpepsi-0}),
%
%
\begin{equation}
\label{eqexpepsi}
\frac\chi k\leq E(\psi_k) \leq\frac\chi{k-1},
\end{equation}
which in turn implies that
\begin{eqnarray*}
E\bigl[\phi_k^{(\ell)}\bigr] &=& E[\psi_k]\prod
_{i=k+1}^{\ell} \bigl(1-E[\psi_i
] \bigr) \leq\frac{\chi}{k-1}\prod_{i=k+1}^\ell
\biggl(1-\frac\chi i \biggr)
\\
&\leq& \frac{\chi}{k-1}\exp\Biggl(-\chi\sum_{i=k+1}^\ell
\frac1i \Biggr) \leq\frac{\chi}{k-1}\exp\biggl(-\chi\log\biggl(
\frac{\ell
+1}{k+1} \biggr) \biggr)
\\
&=& \frac{\chi}{k-1} \biggl(\frac{k+1}{\ell+1} \biggr)^\chi,
\end{eqnarray*}
implying that
%
%
\begin{eqnarray}
\label{eqexp-of-d-up} E\bigl[d_n(k)\bigr] &\leq& m+m\chi
\frac{(k+1)^\chi}{k-1} \sum_{\ell=k}^{n-1} \biggl(
\frac1{\ell+1} \biggr)^\chi\nonumber\\
&\leq& m+m\chi\frac{(k+1)^\chi}{k-1} \int
_{k-1}^{n-1} dx \biggl(\frac1{x+1}
\biggr)^\chi
\nonumber
\\
&=& m+m\frac{\chi}{1-\chi}\frac{(k+1)^\chi}{k-1} \bigl(n^{1-\chi
}-k^{1-\chi}
\bigr)
\\
&\leq& m+ m\frac{\chi}{1-\chi}\frac{k+1}{k-1} \biggl( \biggl(\frac{n}{k}
\biggr)^{1-\chi} -1 \biggr)
\nonumber
\\
&\leq& m+ m\frac{\chi}{1-\chi} \biggl( \biggl(\frac{n}{k}
\biggr)^{1-\chi} -1 \biggr) \biggl(1+\frac4k \biggr).
\nonumber
\end{eqnarray}
On the other hand, again by (\ref{eqexpepsi}),
\begin{eqnarray*}
E\bigl[\phi_k^{(\ell)}\bigr] &\geq& \frac{\chi}{k}\prod
_{i=k+1}^\ell\biggl(1-\frac\chi{ i-1}
\biggr) \geq\frac{\chi}{k}\prod_{i=k+1}^\ell
\biggl(1-\frac1{ i-1} \biggr)^\chi
\\
&=& \frac{\chi}{k}\prod_{i=k+1}^\ell
\biggl(\frac{i-2}{ i-1} \biggr)^\chi=\frac{\chi}{k} \biggl(
\frac{k-1}{ \ell-1} \biggr)^\chi
\end{eqnarray*}
implying that
%
%
\begin{eqnarray}
\label{eqexp-of-d-low} E\bigl[d_n(k)\bigr] &\geq& m+m\chi
\frac{(k-1)^\chi}{k} \sum_{\ell=k}^{n-1} \biggl(
\frac1{\ell-1} \biggr)^\chi\nonumber\\
&\geq& m+m\chi\frac{(k-1)^\chi}{k} \int
_{k}^{n} dx \biggl(\frac1{x-1}
\biggr)^\chi
\nonumber
\\
&=& m+m\frac{\chi}{1-\chi}\frac{(k-1)^\chi}{k} \bigl((n-1)^{1-\chi
}-(k-1)^{1-\chi}
\bigr)
\\
&=& m+m\frac{\chi}{1-\chi}\frac{k-1}{k} \biggl( \biggl(\frac{n-1}{k-1}
\biggr)^{1-\chi}-1 \biggr)
\nonumber
\\
&\geq& m+ m\frac{\chi}{1-\chi} \biggl( \biggl(\frac{n}{k}
\biggr)^{1-\chi} -1 \biggr) \biggl(1-\frac1k \biggr).
\nonumber
\end{eqnarray}

\section{Applications}

\subsection{Degree distribution of an early vertex}
\label{seclim-dn}

In this section, we will show that for $n\gg k\gg1$, $d_n(k)$
grows like $ (\frac nk )^{1-\chi} = (\frac
nk )^{\psi/(\psi+1)}$. To give the precise statement, we
need some definition. To this end, let us consider the random
variables
\[
M_k^{({\ell})}=\prod_{j=k+1}^{\ell}
\frac{1-\psi_j}{1-E[\psi_j]}.
\]
The bounds (\ref{eqexpepsi-0}) and (\ref{beta-moments}) imply
that the second moment of $M_k^{({\ell})}$ is bounded uniformly
in ${\ell}$, so by the martingale convergence theorem,
$M_k^{({\ell})}$ converges both a.s. and in $L^2$.
Since $1-E[\psi_j]= (\frac{j-1}j )^\chi+O(j^{-2})$,
this also implies that the limit
%
%
\begin{equation}
\label{Fk-def} F_k=\lim_{{\ell}\to\infty} \prod
_{j=k+1}^{\ell} (1-\psi_j) \biggl(\frac
j{j-1} \biggr)^\chi=\lim_{{\ell}\to\infty} \biggl(
\frac{\ell}{k} \biggr)^\chi\prod_{j=k+1}^{\ell}
(1-\psi_j)
\end{equation}
exists a.s. and in $L^2$. In the following lemma,
$O_P(k^{-1/2})$ stand for a random variable $A$ such that
$Ak^{1/2}$ is bounded in probability.

\begin{Lemma}\label{lem51}
Consider the sequential model for some $\alpha$ and $m$, and
let $F_k$ be as above. Then
%
%
\begin{equation}
\label{dnk-limit} \frac{d_n(k)}{n^{1-\chi}} \to\frac{m}{1-\chi}k^\chi
\psi_k F_k \qquad\mbox{as } n\to\infty,
\end{equation}
both in expectation and in distribution. Furthermore,
\[
F_k>0 \qquad\mbox{a.s. for all $k\geq1$},\qquad \log F_k
=O_P\bigl(k^{-1/2}\bigr)
\]
and
\[
E[F_k]=1+ O
\bigl( k^{-1}\bigr),
\]
implying in particular that
\[
\lim_{n\to\infty}\frac{E[d_n(k)]}{n^{1-\chi}} = \frac{m\chi}{1-\chi}
\frac{1}{k^{1-\chi}}\bigl(1+O\bigl(k^{-1}\bigr)\bigr).
\]
\end{Lemma}

\begin{Remark*}
Note that (\ref{dnk-limit}) holds also for the independent and the
conditional models. The reason is that by the approximating coupling,
the total variation distance between the degree distribution of vertex
number $k$ in the sequential model and that of vertex number $k$ in the
independent (or conditional) model goes to $0$ as $k$ goes to infinity,
and the convergence is uniform in $n$ (the size of the graph).
\end{Remark*}

\begin{pf*}{Proof of Lemma \ref{lem51}}
We first consider the conditional expectation $E[d_n(k)\mid\FF]$,
where, as before, $\FF$ is the $\sigma$-algebra
generated by $\{\psi_h\}_{h=1}^\infty$. Fix $\varepsilon$, and
let $K$ be such that for $\ell\geq K$,
\[
\Biggl\|F_k- \biggl(\frac{\ell}{k} \biggr)^\chi\prod
_{j=k+1}^{\ell} (1-\psi_j)
\Biggr\|_2\leq\varepsilon.
\]
Bounding
\[
\Biggl\|E\bigl[d_n(k)\mid\FF\bigr]- \sum_{\ell={K}}^{n-1}m
\phi_k^{(\ell)} \Biggr\|_2\leq mK
\]
we then approximate
\begin{eqnarray*}
\sum_{\ell={K}}^{n-1}m\phi_k^{(\ell)}
&=& m\psi_k \sum_{\ell={K}}^{n-1}\prod
_{j=k+1}^{\ell}(1-\psi_j)= m
\psi_k \sum_{\ell={K}}^{n-1} \biggl(
\frac k\ell\biggr)^\chi\bigl(F_k +O(\varepsilon) \bigr)
\\
&=& n^{1-\chi} \biggl(\frac{m}{1-\chi}k^\chi\psi_k
F_k+O(\varepsilon) \biggr),
\end{eqnarray*}
where the errors $O(\varepsilon)$ stand for errors in $L^2$. We
thus have show that as $n\to\infty$,
\[
\frac1{n^{1-\chi}}E\bigl[d_n(k)\mid\FF\bigr]\to
\frac{m}{1-\chi}k^\chi\psi_k F_k
\]
in $L^2$. Taking expectations on both sides, we obtain that
(\ref{dnk-limit}) holds in expectation.

To prove convergence in distribution, it is clearly enough to show that
$E[d_n(k)\mid\FF] - d_n(k)\to0$ in probability. But this follows by
an easy second moment estimate and the observation that
\[
E\bigl[d_n(k)^2\mid\FF\bigr]\leq E\bigl[d_n(k)\mid
\FF\bigr]^2 +E\bigl[d_n(k)\mid\FF\bigr].
\]

Next we observe that the bounds established in
Section \ref{secpolest} imply that there is a constant
$C<\infty$ such that for $k\geq2$,
\[
\bigl|\log M_k^{(\ell)}\bigr|\leq\varepsilon+ \frac Ck
\]
with probability at least $1-\frac{C}{\varepsilon^2k}$. Since
these bounds are uniform in $\ell$, they carry over to the
limit, and imply both that a.s. $F_k>0$ for all fixed $k\geq
2$, and that $\log F_k=O_P(k^{-1/2})$ as $k\to\infty$. To prove
that a.s. $F_1>0$, we note that $F_1/F_2$ is proportional to
$1-\psi_2$. The bound $E[F_k]=1+O(k^{-1})$ finally follows from
the fact that $E[M_k^{(\ell)}]=1$ and the observation that
$1-E[\psi_j]= (\frac{j-1}j )^\chi+O(j^{-2})$.
\end{pf*}

\subsection{Degree distribution}

By Theorem \ref{thmmain} and
Corollary \ref{corlimit}, the limiting degree distribution of
the preferential attachment graph $G_n$ is exactly the degree
distribution of the
root of the P\'olya-point graph. As we will see, this allows us to
explicitly calculate the limiting degree distribution of the
preferential attachment graph. In a similar way, it also allows us to calculate
the limiting degree distribution of a vertex chosen at random from the
vertices that
receive an edge from a uniformly random vertex $v_0$ in $G_n$.
We summarize the results in the following lemma.

\begin{Lemma}\label{lemma5.2}
Let $v_0$ be a uniformly chosen vertex in $G_n$, let $D$
be the degree of $v_0$ and let $D'$ be the degree of a vertex
chosen uniformly at random from the $m$ vertices which received
an edge from $v_0$. In the limit $n\to\infty$, the distribution
of $D$ and $D'$ for all three versions of the preferential
attachment graph converge to
\[
\prob(D=m+k) = \frac{\psi+1}\psi\frac{\Gamma(a+1/\psi+1)}{\Gamma
({a})} \frac{\Gamma(k+{a})}{\Gamma(a+1/\psi+k+2)}
\]
and
\[
\prob\bigl(D'=m+1+k\bigr) =\frac{\psi+1}{\psi^2}\frac{\Gamma(a+
1/\psi+1)}{\Gamma({a+1})}
\frac{(k+1)\Gamma(k+{a+1})}{\Gamma(a+1/\psi+k+3)},
\]
where ${a}=m+2mu$. As $k\to\infty$, this gives
\[
\prob(D=m+k)=Ck^{-2-1/\psi}\bigl(1+O\bigl(k^{-1}\bigr)\bigr)
\]
and
\[
\prob\bigl(D'=m+1+k\bigr)=\tilde Ck^{-1-1/\psi}\bigl(1+O
\bigl(k^{-1}\bigr)\bigr)
\]
for some constants $C$ and $\tilde C$ depending on
$m$ and $\alpha$.
\end{Lemma}

Note that for $\alpha=0$, the statements of the lemma reduce to
\[
\prob(D=m+k)=\frac{2m(m+1)}{{(m+k)(m+k+1)(m+k+2)}}
\]
and
\[
\prob\bigl(D'=m+1+k\bigr)=\frac{2(m+1)(k+1)}{{(m+k+1)(m+k+2)(m+k+3)}}.
\]

\begin{pf*}{Proof of Lemma \ref{lemma5.2}}
First we condition on the position $x_0$ of the root of the P\'olya graph.
Let $D$ be the degree of the
root. $D$ conditioned on $x_0$ is $m$ plus a Poisson variable
with parameter
\[
\frac{\gamma}{x_0^\psi}\int_{x_0}^1{ \psi}
x^{\psi-1}\,dx 
=\gamma\frac{1-x_0^\psi}{x_0^\psi},
\]
where $\gamma$ is a Gamma variable with parameters $a= m+2mu$
and $1$.
Let
\[
\kappa=\kappa(x_0)=\frac{1-x_0^\psi}
{x_0^\psi}.
\]
Then
%
%
\begin{eqnarray}
\label{eqkappaalpha} \prob(D=m+k\mid x_0) &=& \frac{\Gamma(k+{a})}{k!\Gamma
({a})}
\frac{\kappa^k}{(\kappa+1)^{k+{a}}} \nonumber\\
&=& \frac{\Gamma(k+{a})}{k!\Gamma
({a})}\frac{(1-x_0^\psi
)^k}{x_0^{k\psi}}\bigl(x_0^\psi
\bigr)^{k+a}
\\
&=& \frac{\Gamma(k+{a})}{k!\Gamma({a})}\bigl(1-x_0^\psi
\bigr)^kx_0^{a\psi}
\nonumber
\end{eqnarray}
and
\begin{eqnarray*}
\prob(D=m+k) &=&{(\psi+1)}\int_0^1
\prob(D=m+k\mid x_0=x)x^\psi \,dx
\\
&=&(\psi+1)\frac{\Gamma(k+{a})}{k!\Gamma({a})} \int_0^1 {
\bigl(1-x^\psi\bigr)^k x^{{(a+1)}\psi}} \,dx
\\
&=&\frac{\psi+1}\psi\frac{\Gamma(k+{a})}{k!\Gamma({a})} \int_0^1
( {1-y} )^k y^{a+1/\psi} \,dy
\\
&=&\frac{\psi+1}\psi\frac{\Gamma(k+{a})}{\Gamma({a})} \prod_{i=1}^{k+1}
\frac{1}{a+1/\psi+i}
\\
&=&\frac{\psi+1}\psi\frac{\Gamma(k+{a})}{\Gamma({a})} \frac{\Gamma(a+
1/\psi+1)}{\Gamma(a+1/\psi+k+2)}.
\end{eqnarray*}

To calculate the distribution of $D'$, we chose $y_0$ uniformly
at random from $[0,x_0]$. Conditioned on $y_0$, the limiting
degree
$D'$ is equal to $m+1$ plus a Poisson variable with
parameter
\[
\frac{\gamma'}{y_0^\psi}\int_{y_0}^1{ \psi}
x^{\psi-1}\,dx =\gamma'\frac{1-y_0^\psi}{y_0^\psi},
\]
where $\gamma'$ is a Gamma variable with parameters $a+1$ and
$1$. Continuing as before, this gives
%
%
\begin{equation}
\label{eqkappaalpha-} \prob\bigl(D'=m+1+k\mid y_0=y
\bigr) = \frac{\Gamma(k+{a+1})}{k!\Gamma({a+1})}\bigl(1-y^\psi\bigr
)^ky^{(a+1)\psi}
\end{equation}
%
%
\begin{eqnarray*}
\prob\bigl(D'=m+1+k\bigr) &=&{(\psi+1)}\int_0^1dx_0x_0^\psi
\frac1{x_0}\int_0^{x_0} dy\prob
\bigl(D'=m+1+k\mid y_0=y\bigr)
\\
&=&(\psi+1)\frac{\Gamma(k+{a+1})}{k!\Gamma({a+1})} \int_0^1dxx^{\psi-1}
\int_0^{x} dy \bigl(1-y^\psi
\bigr)^ky^{(a+1)\psi}
\\
&=&\frac{\psi+1}{\psi^2}\frac{\Gamma(k+{a+1})}{k!\Gamma({a+1})} \int_0^1du
\int_0^{u} dv (1-v)^kv^{a+1/\psi}.
\end{eqnarray*}
Exchanging the integral over $u$ and $v$ we obtain
\begin{eqnarray*}
\prob\bigl(D'=m+1+k\bigr) &=&\frac{\psi+1}{\psi^2}\frac{\Gamma
(k+{a+1})}{k!\Gamma({a+1})}
\int_0^{1} dv (1-v)^kv^{a+1/\psi}
\int_v^1du
\\
&=&\frac{\psi+1}{\psi^2}\frac{\Gamma(k+{a+1})}{k!\Gamma({a+1})} \int_0^{1}
dv (1-v)^{k+1}v^{a+1/\psi}
\\
&=&\frac{\psi+1}{\psi^2}\frac{(k+1)\Gamma(k+{a+1})}{\Gamma({a+1})} \frac
{\Gamma(a+1/\psi+1)}{\Gamma(a+1/\psi+k+3)}.
\end{eqnarray*}

The asymptotic behavior as $k\to\infty$ follows from the well-known
asymptotic behavior of the Gamma function.
\end{pf*}

\subsection{Joint degree distributions}

We can use the same calculation in order to determine the joint
distribution of the degree of the root of the preferential attachment
graph with a
vertex chosen uniformly among the $m$ vertices that
receive an edge from the root.

\begin{Lemma}\label{lemma5.3}
Let $v_0$ be a uniformly chosen vertex in $G_n$, let $D$ be the
degree of $v_0$ and let $D'$ be the degree of a vertex chosen
uniformly at random from the $m$ vertices which received an
edge from $v_0$. In the limit $n\to\infty$, the joint
distribution of $D$ and $D'$ for all three versions of the
preferential attachment graph converges to
\begin{eqnarray*}
&&
\prob\bigl(D'=m+1+k,D=m+j\bigr) \\
&&\qquad= \frac{\psi+1}{\psi^2}
\frac{\Gamma(k+{a+1})}{k!\Gamma({a+1})}\frac
{\Gamma(j+{a})}{j!\Gamma({a})} \int_0^{1}
{dv} (1-v)^kv^{a+1/\psi}\int_v^1du
(1-u)^ju^{a},
\end{eqnarray*}
where ${a}=m+2mu$. As $k\to\infty$ while $j$ is fixed, this
gives
\[
\prob\bigl(D'=m+1+k\mid D=m+j\bigr)=C_jk^{-1-1/\psi}
\biggl(1+O \biggl(\frac1k \biggr) \biggr),
\]
where $C_j$ is a constant depending on $j$, $m$ and $\alpha$,
while for $k$ fixed and $j\to\infty$, we have
\[
\prob\bigl(D=m+j\mid D'=m+1+k\bigr)= \tilde C_k
j^{-a-3-1/\psi} \biggl(1+O\biggl( \biggl(\frac1j \biggr) \biggr)
\biggr),
\]
where $\tilde C_k$ is a constant depending on $k$, $m$ and
$\alpha$.
\end{Lemma}

Note that the conditioning on $D$ does not change the
power law for the degree distribution of $D'$, while the
conditioning on $D'$ leads to a much faster falloff for the
degree distribution of $D$. Intuitively, this can be explained
by the fact that earlier vertices tend to have higher degree.
Conditioning on the degree $D'$ to be a fixed number therefore
makes it more likely that at least one of the $m$ vertices
receiving an edge from $v_0$ was born late, which in turn makes
it more likely that $v_0$ was born late. This in turn makes it
much less likely that the root $v_0$ has very high degrees,
leading to a faster decay at infinity. This effect does not
happen for the distribution of $D'$ conditioned on $D$, since
the vertices receiving edges from the root are born
\textit{before} the root. Note the fact that the exponent
of the power law of the distribution of $D$ conditioned on $D'$ depends
(through $a$) on
$m$.
Heuristically, this seemingly surprising result follows from the fact that
the distribution of the degree of the vertex at time $k$ is (in the
limit) a discretized Gamma distribution with parameter $a$ (i.e., the
probability of being equal $k$ is proportional to $e^{-k/\lambda}\cdot k^a$.
$\lambda$ here is basically an appropriate power of $n/k$). Note that
with this distribution, when $\lambda$ is relatively large the
probability of the degree being small is approximately $\lambda^{-a}$.
This means that when $D^\prime$ is small, the probability that $k$ is
small (i.e., $n/k$ is large) is as small as $(n/k)^a$. But for $D$ to
be big, $k$ needs to be small (up to an exponential tail).
This is the intuitive explanation for the parameter $a$ comes into
the exponent of
the joint distribution.


\begin{pf*}{Proof of Lemma \ref{lemma5.3}}
Let $x_0$ be the location of the root in the P\'olya-point graph, and let
$y_0$ be the location of a vertex chosen uniformly at random
from the $m$ vertices of type ${ L}$ connected to the root. Then
\begin{eqnarray*}
&&\prob\bigl(D'=k+m+1,D=j+m\bigr)
\\
&&\qquad =(\psi+1)\int_0^1dx x^{\psi-1} \int
_0^x {dy}\prob\bigl(D'=k+m\mid y_0=y
\bigr) \\
&&\qquad\quad\hspace*{49pt}{}\times\prob(D=j+m\mid x_0=x).
\end{eqnarray*}
Using (\ref{eqkappaalpha}) and (\ref{eqkappaalpha-}),
we can write this explicitly a
\begin{eqnarray*}
&&\prob\bigl(D'=k+m+1,D=j+m\bigr)
\\
&&\qquad=(\psi+1)\frac{\Gamma({k+a +1})}{k!\Gamma({ a+1})}\frac{\Gamma
(j+{a})}{j!\Gamma({a})} \int_0^1dx
\bigl(1-x^\psi\bigr)^jx^{(a+1)\psi-1} \\
&&\qquad\quad\hspace*{0pt}{}\times\int
_0^x {dy} \bigl(1-y^\psi
\bigr)^ky^{(a{ +1})\psi}
\\
&&\qquad=\frac{\psi+1}{\psi^2}\frac{\Gamma(k+{a+1})}{k!\Gamma
({a+1})}\frac{\Gamma(j+{a})}{j!\Gamma({a})} \int
_0^1du (1-u)^ju^{a}
\int_0^{u} {dv} (1-v)^kv^{a+1/\psi}
\\
&&\qquad=\frac{\psi+1}{\psi^2}\frac{\Gamma(k+{a+1})}{k!\Gamma
({a+1})}\frac{\Gamma(j+{a})}{j!\Gamma({a})} \int
_0^{1} {dv} (1-v)^kv^{a+1/\psi}
\int_v^1du (1-u)^ju^{a}.
\end{eqnarray*}
We want to approximate the double integral by a product of
integrals. Clearly
\begin{eqnarray*}
&&
\int_0^{1} {dv} (1-v)^kv^{a+1/\psi}
\int_v^1du (1-u)^ju^{a}\\
&&\qquad\leq \int_0^{1} {dv} (1-v)^kv^{a+1/\psi}
\int_0^1du (1-u)^ju^{a}
\\
&&\qquad= k! j!\frac{\Gamma(a+1/\psi+1)}{\Gamma(a+1/\psi+k+2)} \frac
{\Gamma(a+1)}{\Gamma(a+j+2)}:=Z.
\end{eqnarray*}
On the other hand,
\begin{eqnarray*}
&&
\int_0^{1} {dv} (1-v)^kv^{a+1/\psi}
\int_0^vdu (1-u)^ju^{a}\\
&&\qquad\leq \int_0^{1} {dv} (1-v)^kv^{a+1/\psi}
\int_0^vdu u^{a}
\\
&&\qquad= \frac1{a+1} k! \frac{\Gamma(2a+1/\psi+2)}{\Gamma(2a+1/\psi
+k+3)}
\\
&&\qquad=\frac{\Gamma(2a+1/\psi+2)}{\Gamma(a+1/\psi
+1)\Gamma(a+2)} \frac{\Gamma(a+1/\psi+k+2)}{\Gamma(2a+1/\psi
+k+3)} \frac{\Gamma(a+j+2)}{ j!}Z
\\
&&\qquad= O \biggl( \biggl(\frac jk \biggr)^{a+1} \biggr) Z,
\end{eqnarray*}
implying that
%
%
\begin{eqnarray}
\label{DcondD1}\qquad
&&\prob\bigl(D'=k+m+1\mid D=j+m\bigr)\nonumber\\[-8pt]\\[-8pt]
&&\qquad=
\frac{1}{\psi} \frac{\Gamma(k+{a+1})}{\Gamma(a+1/\psi+k+2)} \frac
{\Gamma(a+1/\psi+j+2)}{\Gamma(a+j+2)} \biggl(1+O \biggl(
\biggl(\frac jk \biggr)^{a+1} \biggr) \biggr).\nonumber
\end{eqnarray}

A similar calculation gives
\begin{eqnarray*}
&&
\int_0^1du (1-u)^ju^{a}
\int_0^{u} {dv} (1-v)^kv^{a+1/\psi}\\
&&\qquad
= \frac{j!}{a+1/\psi+1} \frac{\Gamma(2a+1/\psi+2)}{\Gamma
(2a+1/\psi+j+3)} \biggl(1+O \biggl(\frac kj \biggr)
\biggr),
\end{eqnarray*}
which in turn implies that for fixed $k$, as $j$ goes to infinity, we get
%
%
\begin{eqnarray}
\label{DcondD2}
&&\prob\bigl(D=j+m\mid D'=k+m+1\bigr)
\nonumber\\
&&\qquad=\frac{\Gamma(2a+1/\psi+2)}{\Gamma({a})\Gamma(a+1/\psi+2)}
\frac{\Gamma(j+{a})}{\Gamma(2a+1/\psi+j+3)}\\
&&\qquad\quad{}\times \frac{\Gamma(a+
1/\psi+k+3)}{\Gamma(k+2)} \biggl(1+O \biggl(
\frac kj \biggr) \biggr).
\nonumber
\end{eqnarray}
The statements of the lemma describing the decay of
(\ref{DcondD1}) and (\ref{DcondD2}) as (resp.) $k\to\infty$ and
$j\to\infty$ follow from the well-known asymptotics of the
$\Gamma$-function.
\end{pf*}

\subsection{Subgraph frequencies}
\subsubsection{\texorpdfstring{Proof of Lemma \protect\ref{lemsub-G-conv}}
{Proof of Lemma 2.4}}

Let $F$ be a finite graph with vertex set $V(F) = \{v_1,
v_2,\ldots, v_k\}$.
As in Section \ref{secsubgraphfrequency}, let $\inj(F,\mathbf n;G_n)$
be the number of injective maps $\Phi$ from $V(F)$ into $V(G_n)$ that
are homomorphisms and
preserve the degrees. In a similar way, given two rooted graphs $(F,v)$
and $(G,x)$,
let
$\widehat\inj((F,v), \mathbf n; (G_n,x))$ be the number of injective
maps $\Phi$ from
$V(F)$ into $V(G_n)$ that are homomorphisms,
preserve the degrees and map $v$ into $x$. Then $\inj(F, \mathbf n;G_n)$
can be reexpressed as
\[
\inj(F, \mathbf n;G_n) = \sum_{x_1\in V(G_n)}
\widehat{\mathrm{inj}}\bigl((F,v_1), \mathbf n; (G_n,x_1)
\bigr).
\]
Since the diameter of $(F,v_1)$ is at most $k$, its image under a homomorphism
$\Phi$ has diameter at most $k$ as well, which in turn implies that
\[
\frac1n\inj(F, \mathbf n;G_n) = \frac1n \sum
_{x_1\in V(G_n)} \widehat{\mathrm{inj}}\bigl((F,v_1), \mathbf n;
B_{k+1}(G_n,x_1)\bigr).
\]
Given $N$ and $r$, let $\mathcal B_r^{(N)}$ be the set of routed
graphs on $\{1,2,\ldots,N\}$ that have radius $r$ and contain
exactly one\vspace*{1pt} of the representatives from each isomorphism class,
and let $\mathcal B_r=\bigcup_{N=1}^\infty\mathcal B_r^{(N)}$.
Then
\begin{eqnarray*}
\frac1n\inj(F, \mathbf n;G_n) & = & \frac1n \sum
_{x_1\in V(G_n)} \widehat{\mathrm{inj}}\bigl((F,v_1), \mathbf n;
B_{k+1}(G_n,x_1)\bigr)
\\
&=&\sum_{B\in\mathcal B_{{ k}+1}} \widehat{\mathrm{inj}}\bigl((F,v_1),
\mathbf n; B\bigr) {\Pr}_{x_1} \bigl(B_{k+1}(G_n,x_1)
\sim B \bigr),
\end{eqnarray*}
where $\sim$ indicates rooted isomorphisms and the probability
is the probability over rooted balls induced by the random
choice of $x_1\in V(G_n)$.

Since $F$ is connected,
$\widehat{\mathrm{inj}}((F,v_1), \mathbf n; B_{k+1}(G_n,x_1))$
is upper bounded by the constant $C=\max_{1 \leq i
\leq k}(n(i) + d_F(v_i))^{k-1} $. Therefore convergence in the sense
of Benjamini--Schramm implies convergence of the right-hand side,
giving that
%
%
\begin{eqnarray}
\label{hat-F-Calc} \hat t(F,\mathbf n)&:=&\lim_{n\to\infty}
\frac1{|V(G_n)|} \inj(F,\mathbf n;G_n)
\nonumber
\\
&=&\sum_{B\in\mathcal B_{{ k}+1}} \widehat{\mathrm{inj}}\bigl((F,v_1),
\mathbf n; B\bigr) {\Pr} \bigl(B_{k+1}(G,x)\sim B \bigr)
\\
&=& E \bigl[\widehat{\mathrm{inj}}\bigl((F,v_1), \mathbf n; (G,x)\bigr)
\bigr],
\nonumber
\end{eqnarray}
where $E[\cdot]$ denotes expectation over the random choices of the
limit graph
$(G,x)$.

\subsubsection{Convergence in probability}

If $G_n$ is a sequence of random graphs, the subgraph frequencies
${\mathrm{inj}}(F, \mathbf n; G_n,)$ are random numbers as well.
Examining the last proof, one easily sees that the expectation
of these numbers converges if $G_n$ converges in the sense of
Definition \ref{defBS-limit}. For the preferential attachment graph,
this gives
\[
\lim_{n\to\infty}\frac1{|V(G_n)|}E \bigl[ \inj(F,\mathbf
n;G_n) \bigr] = \hat t(F,\mathbf n),
\]
where
%
%
\begin{equation}
\label{t-hat-def} \hat t(F,\mathbf n)= E \bigl[\widehat{\mathrm{inj}}
\bigl((F,v_1), \mathbf n; (T,{ 0})\bigr) \bigr]
\end{equation}
with $(T,0)$ denoting the P\'olya-point graph. It turns out that we
can prove a little more,
namely convergence in probability.

\begin{Lemma} Let $G_n$ be one of the three versions of the
preferential attachment graph
defined in Section \ref{secdef-mod}, let $F$ be a finite connected
graph and let
$\mathbf n:V(F)\to\{0,1,\ldots,\}$.
Then
\[
\frac{1}{n} \inj\bigl((F,\mathbf n);G_n\bigr)\to\hat t(F,\mathbf
n) \qquad\mbox{in probability}.
\]
\end{Lemma}

\begin{pf}
Assume that $x_0$ and $x_0'$ are chosen independently uniformly at
random from
$V(G_n)$. Repeating the proof of Theorem \ref{thmmain}, one easily obtains
that the pair $((G_n,x_0),(G_n,x_0'))$ converges to two independent
copies of
the P\'olya-point graph [more precisely, that the distribution of all
pairs of balls
$(B_r(G_n,x_0),B_r(G_n,x_0'))$ converges to the product distribution of the
corresponding balls in $(T,0)$]. As a consequence, the expectation of
$ [\frac{1}{n} \inj((F,\mathbf n);G_n) ]^2$ converges
$ [\hat t(F,\mathbf n) ]^2$, which in turn implies the claim.
\end{pf}

\subsubsection{Calculation of subgraph frequencies}
\label{secfinite-ball}

In this subsection, we calculate the limiting subgraph frequencies
$\hat t(F,\mathbf n)$ using the expression (\ref{t-hat-def}).
Alternatively, one could use the intermediate expression in (\ref{hat-F-Calc})
and the fact that for each given rooted graph $B$ of radius $k$, we can
calculate the probability that the ball of radius $k$ in the P\'
olya-point graph
$(T,0)$ is isomorphic to $B$. But this gives an expression involving
the countably infinite sum over the balls in $\mathcal B_{{k}+1}$,
while
our calculation below only involves a finite number of terms.


In a preliminary step,
we note that the P\'olya-point graph $(T,0)$ and the
point process $\{x_{\bar a}\}$ can be easily recovered from
the countable graph on $[0,1]$ which is obtained by joining two
points $x,x'\in[0,1]$ by an edge whenever $x=x_{\bar a}$ and
$x'=x_{\bar a'}$ for a pair of neighbors $\bar a,\bar a'$ in
$T$. Identifying the point $x_0$ as the root, we obtain an
infinite, random rooted tree on $[0,1]$
which we will again denote by $T$.

Recalling (\ref{t-hat-def}),
we will want to calculate
the expected number
of
maps $\phi$ from $V(F)$ to $[0,1]$ and are
degree preserving homomorphism
from $(F,\mathbf n)$ into $T$ that map $v_1$ into the root $x_0$.
To this end, we explore the tree structure around the node { $x_0$} in $T$,
in a similar fashion as in
Section \ref{secexpltree}.
Obviously, if $F$ is not a tree, then $\hat t(F,\mathbf n) = 0$.
Otherwise, denote the vertex $v_1 \in V(F)$ as the root and
obtain a rooted tree in which the set of children of every node is
uniquely defined.

A mapping $\phi$ from vertices $v_1, v_2,\ldots, v_k$ to
points $x_1, x_2,\ldots, x_k$ on the interval $[0,1]$ defines a
natural total order $\theta$
on $V(F)$. We say a
mapping is consistent with total order $\theta$ if and only if
for every $i$ and $j$, $\theta(v_i) < \theta(v_j)$ implies $x_i < x_j$.

Given the positions $x_1, x_2,\ldots, x_k$ (or equivalently the
ordering $\theta$),
we can divide the children of every node
$v_i$ to two sets $L(v_i)$ and $R(v_i)$, depending on whether their
corresponding
points on the interval are to the left or right of $x_i$, respectively.
With a slight abuse of notation, define\vspace*{1pt} $L = \bigcup_{1 \leq i \leq k} L(v_i)$
and ${ R} = \bigcup_{1 \leq i \leq k} { R}(v_i)$.
Note that $\{v_2,\ldots,v_k\}$ is the disjoint union of $L$ and $R$.
Since we require that the degrees are preserved, the degree of a node
$x_i$ in $T$
is $d_F(v_i)+n_i$. For the root $x_1=x_0$ this gives $d_F(v_1)+n_1$ children,
$m$ to the left, and $n'_1+|R(v_1)|=d_F(v_1)+n_1-m$ to its right.
If $v_i \in L$, its parent appears on its right.
Therefore, of $n(v_i)$ remaining neighbors of $x_i$ that are not mapped
to any vertex in
$F$, $n'(v_i) = d_F(v_i) + n(v_i) - (m + |R(v_i)| + 1)$
should appear to its right-hand side. For $v_i \in R$, $n'(v_i) =
d_F(v_i) + n(v_i) - (m + |R(v_i)|)$.

Using the above notation, we can finally write the probability density function
$p(F, \mathbf n, x)$ for a mapping from $V(F)$ to $x = (x_1, x_2,\ldots,
x_k)$ to be homomorphic
and degree preserving. Conditioned on $\gamma(x_i)=\gamma_i$, it can
be written as
%
%
\begin{eqnarray}\qquad
&&
p(F, \mathbf n, x, \gamma)\nonumber\\[-8pt]\\[-8pt]
&&\qquad={ (\psi+1)x_1^{\psi}} 
\prod_{v_i\in V}{ \biggl( \frac{\exp(-H_i)H_i^{n'(i)}}{n'(i)!} \prod
_{v_j \in L(v_i)} x_i^{-1} \prod
_{v_j \in R(v_i)} \gamma_i\frac{\psi x_j^{\psi-1}}{x_i^\psi} \biggr)
},\nonumber
\end{eqnarray}
where
\[
H_i=\gamma_i\frac{1-x_i^{\psi}}{x_i^\psi}.
\]
The two inner product terms in the above equations are derived using
the description of the
P\'olya-point in Section \ref{secpolyapointdef}. The first term captures
the probability that the remaining degree of $x_i$ is the desired value $n'(i)$.
Indeed, recalling that the children $x>x_i$ of a vertex $x_i$
are given by a Poison process with density
$ \gamma_i\frac{\psi x^{\psi-1}}{x_i^\psi}$ on $[x_i,1]$, we see that
$n_i'$ is a Poisson random variable
with rate
\[
\gamma_i\int_{x_i}^1
\frac{\psi x^{\psi-1}}{x_i^\psi}\,dx=H_i,
\]
giving the first term in the product above.

Also, $\gamma_i$ is a Gamma variable with parameters $\alpha{
(i)}$ and 1, where $\alpha_i$ depends on whether we discover $v_i$
from right or left.
\[
\alpha(i)=\cases{m+2mu+1, &\quad if $v_i \in L$,
\cr
m+2mu, &\quad if
$v_i \in R$.}
\]

Similarly, $\alpha(1)=m+2mu$. Let $C(\theta)$ be the simplex
containing all points $x = (x_1, x_2,\ldots, x_k)$ consistent
with an ordering $\theta$.
Setting
\[
\hat t(F, \mathbf n, \theta) = \int_{C(\theta) \times(0,\infty)^{k}} \prod
_{i=1}^k \frac{e^{-\gamma_i}\gamma_k^{\alpha_k-1}}{\Gamma
(\alpha_i)} p(F,\mathbf n, x,
\gamma)\,dx_1 \cdots dx_k \,d\gamma_1 \cdots
d\gamma_k,
\]
%
$t(F, \mathbf n)$ can now
be computed by summing
$t(F, \mathbf n, \theta)$
over the
$k!$ choices of $\theta$.

\section*{Acknowledgment}

We thank an anonymous referee for helping us in improving the
presentation of the paper.
The research was performed while N. Berger and A. Saberi were visiting
Microsoft Research.



\printaddresses

\end{document}